\documentclass[times,final,11pt]{elsarticle}
\usepackage{amsmath}
\usepackage{amsfonts}
\usepackage{fullpage}
\usepackage[dvipsnames]{xcolor}
\usepackage{subfig} 
\usepackage{graphicx}
\usepackage{subcaption}
\usepackage{caption}
\usepackage{mathtools}
\usepackage{multirow} 
\usepackage[T1]{fontenc}
\usepackage[utf8]{inputenc}

\usepackage{algorithm}
\usepackage{algpseudocode}
\usepackage[colorlinks=true, linkcolor=blue, citecolor=blue, urlcolor=blue]{hyperref}
\newcommand{\bx}{\boldsymbol{x}}

\begin{document}
%

\title{Parallel-in-Time Solution of Allen–Cahn Equations by Integrating Operator Learning into the Parareal Method}

\author[1]{Yuwei Geng}
\ead{ygeng@email.sc.edu}
\author[2]{Junqi Yin}
\ead{yinj@ornl.gov}
\author[3]{Eric C. Cyr}
\ead{eccyr@sandia.gov}
\author[4]{{Guannan Zhang}}
\ead{zhangg@ornl.gov}
\author[1]{{Lili Ju}\texorpdfstring{\corref{cor1}}{}}
\ead{ju@math.sc.edu}

\cortext[cor1]{Corresponding author}
\address[1]{Department of Mathematics, University of South Carolina, Columbia, SC 29208, USA}
\address[2]{National Center for Computational Sciences, Oak Ridge National Laboratory, TN 37831, USA}
\address[3]{Computational Mathematics Department, Sandia National Laboratories, NM 87123, USA}
\address[4]{Computer Science and Mathematics Division, Oak Ridge National Laboratory, TN 37831, USA}

\begin{frontmatter}
\begin{abstract}
While recent advances in deep learning have shown promising efficiency gains in solving time-dependent partial differential equations (PDEs), matching the accuracy of conventional numerical solvers still remains a challenge. One strategy to improve the accuracy of deep learning-based solutions for time-dependent PDEs is to use the learned solution as the coarse propagator in the Parareal method and a traditional numerical method as the fine solver. However, successful integration of deep learning into the Parareal method requires consistency between the coarse and fine solvers, particularly for PDEs exhibiting rapid changes such as sharp transitions. To ensure such consistency,
we propose to use the convolutional neural networks (CNNs) to learn the fully discrete time-stepping operator defined by the traditional numerical scheme used as the fine solver. We demonstrate the effectiveness of the proposed method in solving the classical and mass-conservative Allen-Cahn (AC) equations.
Through iterative updates in the Parareal algorithm, our approach achieves a significant computational speedup compared to traditional fine solvers while converging to high-accuracy solutions. Our results highlight that the proposed Parareal algorithm effectively accelerates simulations, particularly when implemented on multiple GPUs, and converges to the desired accuracy in only a few iterations. Another advantage of our method is that the CNNs model is trained on trajectories-based on random initial conditions, such that the trained model can be used to solve the AC equations with various initial conditions without re-training.
This work demonstrates the potential of integrating neural network methods into the parallel-in-time frameworks for efficient and accurate simulations of time-dependent PDEs.
\end{abstract}


\end{frontmatter}

\section{Introduction}
The Allen-Cahn (AC) equation was originally proposed in \cite{allen1979microscopic} as a phenomenological model to describe antiphase domain coarsening in binary alloys. Formulated as a time-dependent partial differential equation, it defines an initial value problem where the temporal evolution of an order parameter captures phase transitions. Beyond its original context, the AC equation has been widely applied in modeling dynamic processes across various scientific and engineering domains, including crystal growth \cite{wheeler1992phase}, image processing \cite{benevs2004geometrical}, and biological pattern formation \cite{shao2010computational}.
Mathematically, analytical solution to the AC equation are very challenging or even impossible to find. Hence, effective and accurate numerical methods for computing their approximate solutions become highly necessary in order to simulate and understand their dynamics. Typically, those are solved through numerical integration over time steps using, e.g., fully implicit schemes~\cite{feng2003numerical,xu2019stability} or stabilized semi-implicit schemes~\cite{shen2010numerical,tang2016implicit}. However, solving the AC equation requires very fine spatial and temporal discretization, with mesh sizes smaller than the interface parameter and time steps restricted by stability conditions, making direct simulation costly. Although implicit methods relax time step restrictions, they still involve solving large systems, so simulations become particularly expensive when having sharp interfaces between phases. Recently, deep learning has provided a promising alternative in solve the AC equation, e.g., \cite{geng2024deep} to significantly improve the computational efficiency. Nevertheless, deep-learning-based solutions of the AC equations often lack sufficient accuracy, as they struggle to capture the dynamics of sharp interface structures.

The overarching goal of this work is to integrate deep-learning-based AC solvers with traditional numerical solvers in a parallel-in-time framework, thereby accelerating the solution of both classic and mass-conservative AC equations while maintaining satisfactory accuracy. 
The Parareal method \cite{lions2001resolution} is a parallel-in-time algorithm designed to accelerate the solution of time-dependent problems by leveraging parallel computing resources. The Parareal algorithm iteratively combines two propagators: a coarse propagator that is computationally efficient and provides an approximate solution, and a fine propagator that is more accurate but computationally expensive. Each iteration uses corrections from the fine solver to refine the coarse prediction, gradually converging to a high-fidelity solution. 

There are two emerging directions linking parallel-in-time integration with machine learning. One direction focuses on applying parallel-in-time principles to accelerate the training of deep learning models. A multigrid reduction-in-time \cite{falgout2014parallel} approach is used to enable layer-parallel training of deep residual networks \cite{gunther2020layer} by interpreting layers as a time domain and applying parallel-in-time techniques for efficient parameter updates.  Extensions of this idea in \cite{kirby2020layer}
have demonstrated additional gains on multi-GPU platforms. 
Another strategy \cite{meng2020ppinn} aims to improve the training efficiency of physics-informed neural networks, especially for long time horizons. By breaking the time domain into shorter segments and applying Parareal to connect them, training becomes more scalable. Other work includes adapting neural ODEs to a parallel-in-time framework to improve residual network training in \cite{lorin2020derivation}, and proposing a Parareal-inspired neural network architecture for efficient parallel training of deep models across multiple GPUs in \cite{lee2022Parareal}.

Another direction of research focuses on enhancing the performance of parallel-in-time algorithms using machine learning. For example, neural networks have been explored as coarse propagators in the Parareal algorithm in a series of efforts \cite{yalla2018parallel,agboh2020Parareal,ibrahim2023Parareal,nguyen2023numerical}. 
Despite the initial success of combining Parareal with deep learning for solving PDEs, there are significant challenges in extending this approach to the classic and mass-conservative AC equations. First, it is difficult to ensure consistency between the deep-learning model used as the coarse propagator and the traditional numerical solver used as the fine propagator, particularly for the AC equation which has dynamically evolving sharp interfaces. Failing to ensure such consistency will require more global iterations in the Parareal algorithm, ultimately diminishing the computational benefits of using Parareal.
Second, the deep-learning model should generalize across a wide range of initial conditions beyond those used in training. Otherwise, the need to re-train the model for each new initial condition would offset the computational benefits of using Parareal due to the added training cost.

To address these challenges, we propose training two neural network models, i.e., an Allen-Cahn Neural Network (ACNN) for the classic AC equation and a mass-conservative Allen-Cahn Neural Network (mACNN) for the mass-conservative AC equation, and integrating them into the Parareal framework to accelerate time-to-solution for both types of the AC equations. The ACNN and mACNN models learn from fully discrete time-stepping operators for initial value problems targeting the classic and mass-conservative AC equations. As operator-learning approaches, they generalize across different initial conditions and therefore do not require retraining when the input data changes. Since the models are trained using the fully discrete time-stepping scheme, they naturally incorporate nonlocal terms and remain consistent with the fine solver used in Parareal. 
Additionally, the neural networks can be executed with high computational efficiency, which enables simulations on sufficiently fine temporal meshes to ensure stability while only requiring coarse-mesh snapshots for Parareal iterations.
Moreover, the convolutional neural network architecture enables them to capture sharp interfacial structures that are often difficult for other machine-learning-based solvers. We evaluated the proposed framework on both the classic and conservative AC equations and demonstrated its scalability across multiple compute nodes, achieving performance on up to 128 GPUs.

The rest of the paper is organized as follows. In Section \ref{sec:prework}, we introduce the numerical scheme for the AC equations and review the ACNN and mACNN models. Section \ref{sec:method} discusses the Parareal algorithm and explains how it integrates with our deep-learning models as coarse propagators. In Section \ref{sec:results}, we present comparisons in terms of accuracy and speedup to demonstrate the effectiveness of the Parareal algorithm. Finally, we conclude the paper with a brief discussion on future directions.

\section{Problem setting}\label{sec:prework}
We introduce the classical AC and mass-conservative AC (mAC) models, followed by the numerical schemes used for both equations. We then describe the network architectures for ACNN and mACNN.  
Let $\Omega \subset \mathbb{R}^d$, $d=2,3$ be a bounded domain. The classic AC equation given by 
\begin{equation}
    \partial_t u = \epsilon^2 \Delta u + f(u)
    \label{eq:class_AC_eq}
\end{equation}
and equipped with periodic boundary conditions and a suitable initial condition $u(0,\bx)=u_0(\bx)$, $\epsilon$ is the interfacial parameter, $u(\bx,t)$ is the unknown scalar function and $f=-F'$, where $F(u)$ is an nonlinear potential function.
The commonly used double-well potential function is adopted for this work:
\begin{equation}
    F(u) = \frac{1}{4}(u^2-1)^2, \quad f(u)=u-u^3.
\end{equation}
Based on the maximum bound principle (MBP) \cite{du2021maximum,li2021unconditionally}, the solution of the classic AC equation \eqref{eq:class_AC_eq} always stays inside the interval $[-1,1]$. 

The classic AC equation does not preserve the total mass, defined as 
\begin{equation}
    M(t) = \int_{\Omega} u(\boldsymbol{x},t) d\boldsymbol{x}.
\end{equation}
To overcome this limitation, several mass \textit{conservative} variants of the AC equation have been developed. A prominent example, proposed and analyzed in \cite{rubinstein1992nonlocal}, introduces a nonlocal modification of the classical formulation: 
\begin{equation}\label{eq:mass_AC_eq}
    \partial_t u = \epsilon^2 \Delta u + f(u)-\frac{1}{|\Omega|}\int_{\Omega} f(u(\boldsymbol{y},t)) d\boldsymbol{y},\quad \boldsymbol{x} \in \Omega, t > 0,
\end{equation}
where the term \( \frac{1}{|\Omega|} \int_{\Omega} f(u(\boldsymbol{y},t)) d\boldsymbol{y} \) acts as a Lagrange multiplier to ensure mass conservation. This guarantees that \( M(t) = M(0) \) for all \( t \geq 0 \). Both the classic AC equation \eqref{eq:class_AC_eq} and its conservative counterpart \eqref{eq:mass_AC_eq} are \textit{autonomous} systems under periodic or homogeneous boundary conditions. The MBP principle implies that solutions of the conservative AC equation  \eqref{eq:mass_AC_eq} remain within the range $[-\frac{2}{3}\sqrt{3},\frac{2}{3}\sqrt{3}]$.
From a variational viewpoint, the classic AC equation is the \( L_2 \) gradient flow associated with the energy functional:  
\begin{equation}\label{eq1:energy}
    E(u) = \int_\Omega \left(\frac{\epsilon^2}{2} |\nabla u(\boldsymbol{x},t)|^2 + F(u(\boldsymbol{x},t))\right) d\boldsymbol{x}.
\end{equation}  
The conservative AC equation \eqref{eq:mass_AC_eq}, while enforcing mass conservation, also retains an energy dissipation property with respect to the same functional \eqref{eq1:energy}.

When solving the AC equation with traditional numerical schemes, the spatial discretization requires the computational mesh size to be smaller than the interfacial parameter $\epsilon$ in order to accurately resolve the transition layers and capture the correct dynamics. 
%
In recent years, neural network models have gained increasing popularity for solving PDEs, owing to their flexibility and ability to approximate complex nonlinear dynamics. However, despite these advances, several notable challenges remain. First, many neural network models—such as standard PINNs \cite{raissi2019physics}—require retraining whenever the initial condition changes. Second, the presence of nonlocal terms in the mAC equation, particularly in higher dimensions, poses significant difficulties for conventional neural architectures. Third, capturing sharp interfacial structures remains a persistent challenge \cite{krishnapriyan2021characterizing,xu2025overview}, as neural networks often struggle with irregular or non-smooth solutions. Although a number of modified models have been proposed to address individual aspects of these difficulties \cite{mattey2022novel,zhao2020solving,chen2025learn,wang2022modified,huang2025frequency,liu2024mitigating}, no single framework has successfully resolved all of them simultaneously.

\section{Operator-Learning-Enabled Parareal Algorithm for Allen-Cahn Equations}\label{sec:method}
This section is organized as follows. We begin by presenting the fully discrete schemes introduced in Section~\ref{sec:CNscheme}. Next, Section~\ref{sec:CNN_model} introduces the neural network model, including its architecture, loss function, and training procedure. Finally, in Section~\ref{sec:Parareal_alg}, we employ the ACNN/mACNN as a coarse propagator within the Parareal framework.

\subsection{The fully-discrete schemes for the Allen-Cahn equations}\label{sec:CNscheme}
For notational simplicity, we present our methodology in a 2D setting, while the numerical experiments in Section \ref{sec:results} include both 2D and 3D cases.
Specifically, we consider a two-dimensional square domain $\Omega = [0, L]^2$ and assume periodic boundary conditions for the solution $u$, 
Note that the numerical scheme is employed not only as the fine solver but also as the basis for defining the loss function used to train the neural network model.
 Let us uniformly partition the domain $\Omega$ and obtain the set of grid points $\{(x_i,y_j)\}_{i,j=0}^N$ where $N>0$ is an integer, $x_i=ih$ and  $y_j=jh$ with $h=L/N$. Then we choose a uniform time step size $\Delta t>0$ and set $t_n = n\Delta t$ for $n\geq 0$. 
Correspondingly, we set $U_0^{i,j}= u_0(x_i,y_j)$ and assume $U_n^{i,j}\approx u(x_i,y_j,t_n)$. The periodic boundary condition implies $U_n^{0,j} = U_n^{N,j}$ and $U_n^{j,0} = U_n^{j,N}$ for $j=0,1,\cdots N$, thus the set of unknown can be represented by $U_n=\{U_n^{i,j}\}_{i,j=1}^N$.

We adopt the second-order Crank-Nicolson approximation in time to discrete the model in \eqref{eq:class_AC_eq} and \eqref{eq:mass_AC_eq}. 
To solve AC equations with fully discretization, we can use central finite-difference discretization or a Fourier collocation method as described in \cite{yang2018uniform}.
With the periodic boundary condition, the numerical scheme can be written as:
\begin{equation*}
     \frac{U^{i,j}_{n+1}-U^{i,j}_{n}}{\Delta t} =  \epsilon^2\Delta\frac{U^{i,j}_{n+1}+U^{i,j}_{n}}{2} +\frac{f(U^{i,j}_{n+1})+f(U^{i,j}_{n})}{2} - g(U_{n},U_{n+1}), \label{eq:nm0_1}
\end{equation*}
for $i,j=1,2,\cdots,N$. 
\begin{eqnarray*}
g(U_{n},U_{n+1})=
\begin{cases}
    0, ~~\,\,\,\, \text{ for } ~~\text{class AC equation},  \\
    \frac{1}{N^2}\sum_{i',j'=1}^N \frac{f(U_{n+1}^{i',j'})+f(U_{n}^{i',j'})}{2},~~\text{ for}~~  \text{mass-conservative AC equation}.
\end{cases}
\end{eqnarray*}
This yields 
\begin{equation}\label{eq:fully-discrete}
     U^{i,j}_{n+1} -
    \frac{\Delta t}{2} \epsilon^2\Delta U^{i,j}_{n+1}=U^{i,j}_{n}+\frac{\Delta t}{2}\left(\epsilon^2\Delta {U^{i,j}_{n}} + {f(U^{i,j}_{n})+ {f(U^{i,j}_{n+1})}-g(U_{n},U_{n+1})}\right),
\end{equation}
for $i,j=1,2,\cdots,N$.
It is known that the system \eqref{eq:fully-discrete} is  guaranteed to have a  unique solution $U_{n+1}$ for any given $U_n$ 
under a quite relaxed restriction on the time step size $\Delta t$ and the resulting numerical solution is of second-order accuracy in both time and space. With picard iteration, the numerical scheme can be solve as: let $U_{n+1,0} = U_n$ and for $m=0,1,\cdots$, compute $U_{n+1,(m+1)}$ by solving the linear system
\begin{equation*}
     \left(I-\frac{\Delta t}{2}\epsilon^2\Delta \right) U^{i,j}_{n+1,{(m+1)}}=U^{i,j}_{n}+\frac{\Delta t}{2}\left(\epsilon^2\Delta {U^{i,j}_{n}} + {f(U^{i,j}_{n})+ {f(U^{i,j}_{n+1,(m)})}-g(U_{n},U_{{n+1},(m)})}\right), 
\end{equation*}
for $i,j=1,2,\cdots,N$. and finally set $U_{n+1} = U_{n+1,m^*}$ if $U_{n+1,m^*}$ satisfies certain convergence criteria. 

The system can be rewritten as 
\begin{equation}\label{eq:fine_operator}
U_{n+1}={\mathcal{F}}(t_n,t_{n+1},U_n), \quad n=0,1,\cdots,
\end{equation}
where ${\mathcal{F}}$ denotes the fully discrete operator that maps the input $U_n$ to the output $U_{n+1}$ through the relation \eqref{eq:fully-discrete} with corresponding $\Delta t$. This operator ${\mathcal{F}}$ serves as the fine propagator in the Parareal algorithm. 

\subsection{Training the ACNN and mACNN models}\label{sec:CNN_model}
Building on the numerical scheme, we introduce the machine learning models employed as coarse propagators within the Parareal framework. Recall that the fine propagator $\mathcal{F}$, introduced in Section~\ref{sec:CNscheme}, denotes the nonlinear evolution operator defined by the fully discrete numerical scheme in equation~\eqref{eq:fine_operator}. For a given state $U_n$, the operator $\mathcal{F}$ advances the solution to the next state $U_{n+1}$.  
In contrast, the nonlinear mapping $\mathcal{G}$ defined by the ACNN or mACNN model
\[
\mathcal{G}(U_n;\Theta): U_n \mapsto U_{n+1},
\]  
where $\Theta$ represents the set of trainable parameters. The model $\mathcal{G}$ is trained so that  
\[
\mathcal{G}(U_n;\Theta) \approx \mathcal{F}(t_n,t_{n+1},U_n),
\]  
thereby providing a approximation of the numerical evolution operator. Within the Parareal algorithm, $\mathcal{F}$ serves as the fine propagator, while $\mathcal{G}$ is used as the coarse propagator, and the interaction between these two operators forms the foundation of the proposed framework.

\begin{figure}[ht!]
    \centering
\includegraphics[width=0.95\linewidth]{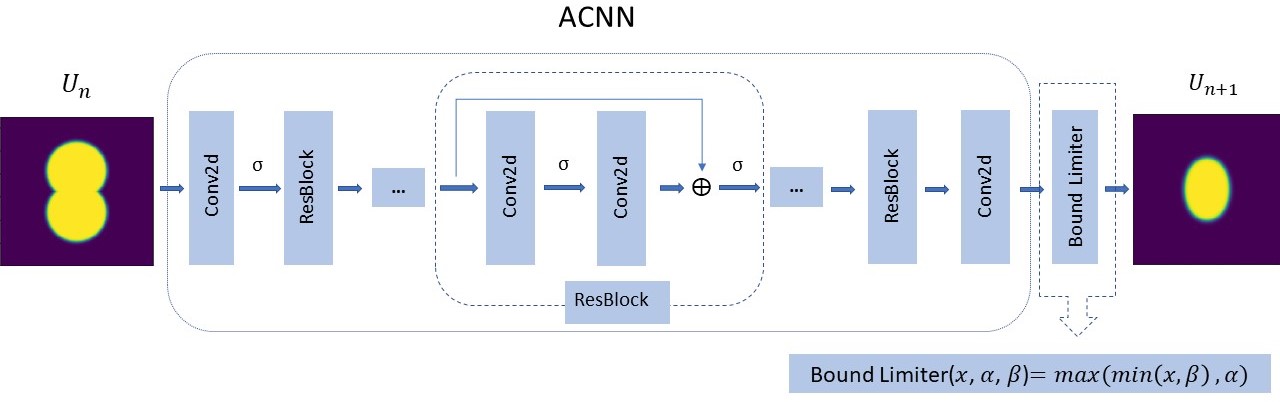}
\includegraphics[width=0.95\linewidth]{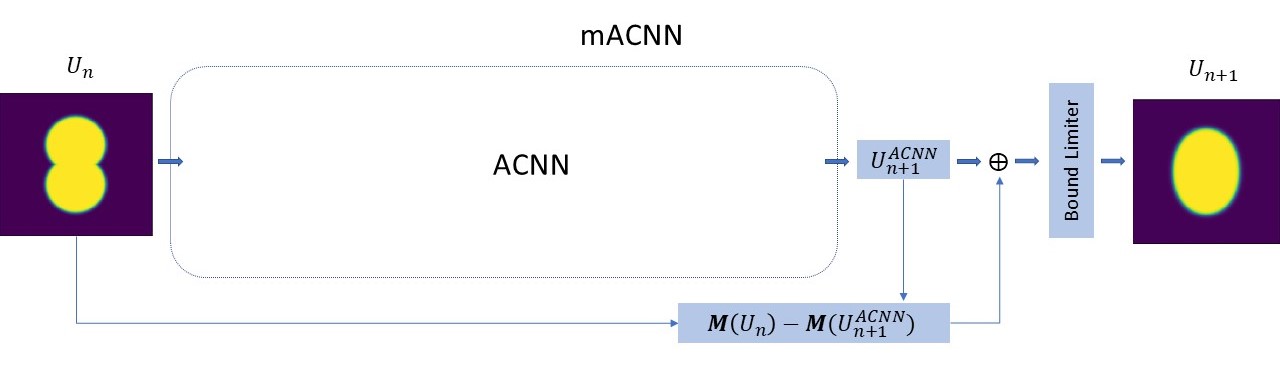}
    \caption{The structure of ACNN (top) and mACNN (bottom) which learns the dynamics of classic AC equation and mass-conservative AC equation respectively.}
    \label{fig:network architecture}
\end{figure}

Figure~\ref{fig:network architecture} illustrates the overall architecture of ACNN (top) and mACNN (bottom).
The input $U_n$ is stored as a four-dimensional tensor of shape $[1,1,N,N]$, corresponding to the format \texttt{[batch, channel, height, width]}, and the output is the predicted state $U_{n+1}$. 
 The network architecture comprises an input layer, several residual blocks (ResBlocks) \cite{targ2016resnet}, an output layer, and a boundary-limiter module. Importantly, both ACNN and mACNN avoid the use pooling or fully connected layers. To respect the periodic boundary conditions inherent to the problem, all convolutional layers employ ``circular'' padding. Each convolutional layer utilizes $3 \times 3$ filters with a padding size of 1. The hyperbolic tangent (tanh) function serves as the activation function after each convolutional layer, except for the output layer. Following the input layer, the tensor shape become $[C,1,N,N]$, which is maintained across all intermediate layers. The finial output tensor $U_{n+1}$ has the shape $[1,1,N,N]$ consistent with the input format. The bound limiter is to enforce the maximum principle (MBP), ensuring that the output remains within a prescribed interval $[\alpha, \beta]$. To solve the mass-conservative AC equation, an additional value is uniformly added to all entries of the network's output to enforce mass conservation before the bound-limited laryer.
Note that in the referenced article, the authors used a more complex architecture with three residual blocks and 16 channels ($C$=16) to achieve high-accuracy simulations. In contrast, our work emphasizes a trade-off between accuracy and computational efficiency, which is particularly important in the context of the Parareal algorithm. For the coarse propagator, it is crucial to maintain reasonable accuracy while ensuring fast inference to maximize parallel efficiency. Therefore, we adopt a lighter network with two residual blocks and 4 channels ($C$=4). Despite its reduced complexity, this model still delivers sufficiently accurate simulations, making it well-suited as a coarse propagator in the Parareal framework.

With the defined of the network architecture, the loss function can be defined as follows. Assuming $R$ different problems are trained simultaneously with batch size $R$, the loss function of the neural network model is defined using the fully discrete scheme~\eqref{eq:fully-discrete} and measured in the $L^2$-norm:  
\begin{equation}\label{eq:loss_ACN} 
\begin{aligned}
\mathcal{L} =& \frac{1}{R N^2}\sum_{r=1}^{R} \sum_{i,j=1}^{N} \left[ U_{n+1}^{i,j,(r)} - U_n^{i,j,(r)} - \frac{\Delta t}{2} \Big( \epsilon^2 \Delta_h (U_{n+1}^{i,j,(r)} + U_n^{i,j,(r)}) + f(U_{n+1}^{i,j,(r)}) + f(U_n^{i,j,(r)})\right.\\
&\qquad\qquad\qquad\left.- g(U_{n+1}^{i,j,(r)}) - g(U_n^{i,j,(r)}) \Big) \right]^2.
\end{aligned}
\end{equation}
Because the loss function is derived directly from the fully discrete numerical scheme, the neural network model is consistent with the fine solver in the sense that it enforces the same discrete evolution law. This consistency ensures that the learned coarse propagator $\mathcal{G}$ remains faithful to the underlying discretization of the AC equations, making it particularly well-suited for integration into the Parareal algorithm, where alignment between fine and coarse solvers is critical for stability and convergence. Moreover, since the formulation involves the fully discrete representations of both $U_n$ and $U_{n+1}$, the model can naturally accommodate nonlocal terms, thereby extending its applicability to the mAC equation.

To train the ACNN and mACNN models following Algorithm~\ref{alg:training}, the training data are generated from multiple problems with different initial conditions, which are randomly partitioned into subsets and evolved in time. The training is carried out iteratively across time steps, with the network parameters updated using the scheme-based loss function~\eqref{eq:loss_ACN}. Consequently, the learned models are not restricted to a single trajectory but are capable of generalizing across arbitrary initial conditions. This property makes the neural network propagator particularly suitable for the Parareal algorithm, where the coarse solver can be applicable to a wide range of solution states.  

We have introduced both the fine propagator $\mathcal{F}$, derived from the fully discrete numerical scheme, and the coarse propagator $\mathcal{G}$, defined by the neural network models ACNN and mACNN, whose loss functions are constructed directly from the fully discrete formulation. In the next section, we combine these two propagators within the Parareal framework and present the resulting time-parallel algorithm in detail.

\begin{algorithm}[H]
\caption{Training strategy of ACNN/mACNN}
\label{alg:training}
\begin{algorithmic}[1]
\State Generate training data from $S$ different problems with distinct initial conditions $\{U^{(r)}_0\}_{r=1}^R$.
\State Randomly partition the initial conditions into $p$ subsets, each containing $q$ samples ($R=pq$).
\For{each subset of initial conditions}
     \State Initialize $U_0$ with $q$ initial condition samples.
    \For{each time step $n=0,1,\dots, T_{\text{train}}/\Delta t$}
        \For{$k=1,\dots,b$} \Comment{inner training loop}
            \State Predict the next state: $U_{n+1} = \mathcal{G}(U_n;\Theta)$.
            \State Update network parameters $\Theta$ by minimizing the loss~\eqref{eq:loss_ACN}.
        \EndFor
        \State After $b$ updates, fix $U_{n+1}$ and use it as input for the next step.
    \EndFor
\EndFor
\State Repeat for multiple epochs until convergence.
\end{algorithmic}
\end{algorithm}

\subsection{Integrating the trained ACNN and mACNN into the Parareal algorithm}\label{sec:Parareal_alg}
\begin{figure}[h!]
    \centering
    \includegraphics[width=0.98\linewidth]{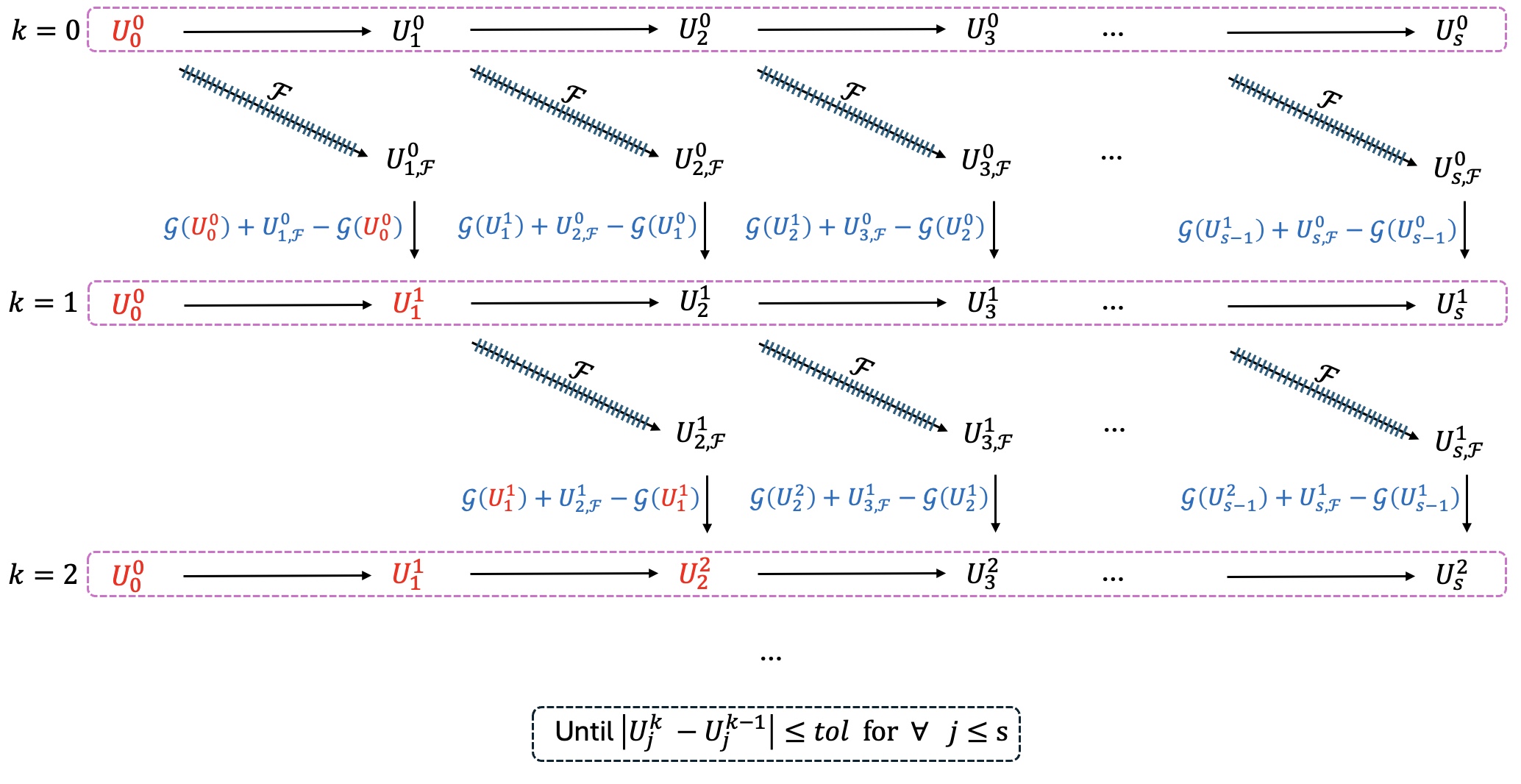}
    \caption{Illustration of the Parareal algorithm. The coarse propagator here is labeled $\mathcal{G}$, whereas the fine propagator is labeled $\mathcal{F}$. The red font means that solutions are already the fine solver simulations.}
    \label{fig:Parareal}
\end{figure}

We now introduce the Parareal algorithm using ACNN and mACNN as the coarse propagator and analyze the speedup. For simplicity, we write $ \mathcal{G}(U_n) = \mathcal{G}(U_n,\Theta) =U_{n+1}$. Figure~\ref{fig:Parareal} illustrates the workflow of the Parareal algorithm.
The algorithm begins with an initial coarse prediction 
\[
\mathbf{U}^0 = [U^0_0, U^0_1, \ldots, U^0_s],
\] 
where $U^0_0 = U_0$ is the prescribed initial condition, the coarse time step size is $\Delta t$ and $U^0_{n+1} =\mathcal{G}(U^0_{n})$ for $n=0,1,\cdots,s-1$. This initial trajectory $\mathbf{U}^0$ is then refined by the fine solver $\mathcal{F}$, which computes high-fidelity solutions in parallel as
\[
U_{n+1,\mathcal{F}}^{0} = \mathcal{F}(t_n, t_{n+1}, U_n^0), \qquad n = 0,1,\ldots,s-1,
\]
using a smaller fine time step $\Delta t^*$. This results in the intermediate states $[U_0^0, U_{1,\mathcal{F}}^{0}, \ldots, U_{s,\mathcal{F}}^{0}]$.  

The updated solutions in the Parareal algorithm are computed iteratively using the formula:
\begin{equation} U_{n+1}^{k} = \mathcal{G}(U^{k}_{n}) + \mathcal{F}(t_n,t_{n+1},U^{k-1}_{n}) - \mathcal{G}(U^{k-1}_{n}), 
\end{equation}
until the sequence of iterates converges. Specifically, we declare convergence at iteration 
 $k$ if the following condition is satisfied $\forall j \leq s$, $|U_{j}^{k} - U_{j}^{k-1}| < tol$
where $tol$ is a user-defined tolerance threshold.

Since the initial condition $U_0^0$ is fixed, after first iteration, we have:
$$
U_1^1 =\mathcal{G}(U^1_{0}) + \mathcal{F}(t_0,t_1,U^0_{0}) - \mathcal{G}(U^0_{0})= \mathcal{F}(t_0,t_1,U^0_{0}).$$
where $U^1_{0} = U^0_{0} = U_0.$
Thus, after the first update, $U_1^1$ exactly matches the solution obtained from the fine solver using time step $\Delta t^*$. Consequentially, after $k$ iterations, the first $k$ time-step solutions 
$[U^0_0,U^1_1, \cdots,U^k_{k}]$ will already match the outputs of the fine solver. Therefore, at most $s$ iterations are required to obtain the full fine-resolution solution over $s$ time steps where iteration times $k\leq s$.

\subsection{Discussion on the computational efficiency}
The time complexity of the Parareal algorithm has been discussed in detail in \cite{minion2011hybrid}. We analyze the computational cost of using a CNN-based coarse solver within the Parareal framework.
Assume the cost of $U_{n+1} = \mathcal{G}(U_n)$ with fixed $\Delta t$ is $t_{nn}$, and the cost using the fine numerical solver $U_{n+1} = \mathcal{F}(t_n,t_{n+1},U_n)$ with fine time step $\Delta t^*$ is $ct_{num}$, where $c = \frac{\Delta t}{\Delta t^*}$. 
Suppose we aim to simulate the solution from $t=0$ to $t=s\Delta t = sc\Delta t^*$. Assuming the time cost per step is constant for both ACNN/mACNN models and numerical solver, the total time cost of a full simulation using only the fine solver is:
$T_{num} =cst_{num}$. 

Assume that the Parareal algorithm terminates after $k$ iterations and the number of GPUs is sufficient to parallelize, then each iteration's parallel cost becomes $ct_{num}$ and the total cost simplifies to:
\begin{align*}
    T(k) &=  ckt_{num}+ \left(s+(s-1)+\cdots+(s-k)\right)t_{nn}\\
    & = ckt_{num}+ \frac{(2s-k)(k+1)}{2}t_{nn}
\end{align*}
Comparing with the baseline cost $T_{num} = cst_{num}$, the Parereal approach is faster if: 
\begin{equation*}
    \frac{(2s-k)(k+1)}{2}t_{nn} \leq (s-k)ct_{num}.
\end{equation*}
The corresponding speedup is: 
\begin{equation*}
    S_k = \frac{T_{num}}{T(k)} = \frac{1}{\frac{k}{s}+ \frac{(2s-k)(k+1)}{2cs} \frac{t_{nn}}{t_{num}}} \leq \text{min}\{\frac{s}{k}, \frac{2cst_{num}}{(2s-k)(k+1)t_{nn}} \}
\end{equation*}
\begin{algorithm}[H]
\caption{The Parareal Algorithm with ACNN/mACNN as the Coarse Solver}
\label{alg:Parareal}
\begin{algorithmic}[1]
\Require Initial condition \( {U}^0 \), number of time steps \( s \), tolerance \( {tol} \)
\State Generate initial guess \( [U^0_0, U^0_1, \dots, U^0_s] \) sequential using the ACNN/mACNN model:
\State \(\quad U^0_{n+1} = \mathcal{G}(U^0_n) \) for \( n = 0, 1, \dots, s-1 \)
\For{$k = 1$ to $s$}
    \State Compute intermediate fine solutions in parallel using fine solver:
    \State \(\quad U^{k-1}_{n+1,\mathcal{F}} = \mathcal{F}(t_n,t_{n+1},U^{k-1}_n) \) for \( n = k-1, k, \dots, s-1 \) \Comment{Using fine time step $\Delta t^*$}
    \For{$n = k-1$ to $s-1$} 
        \State \( U^k_{n+1} = \mathcal{G}(U^{k}_{n}) + U^{k-1}_{n+1,\mathcal{F}}  - \mathcal{G}(U^{k-1}_{n}) \)
        \Comment{Update solution using Parareal correction}
    \EndFor
    \If{$|U^k_j - U^{k-1}_j| < {tol}$ for all $j \leq s$}
        \State Convergence achieved, exit loop.
    \EndIf
\EndFor
\Ensure Converged solution \( [U^0_0, U^1_1,\cdots,U^k_k, \dots, U^k_s] \)
\end{algorithmic}
\end{algorithm}

The efficiency of the proposed Parareal framework with ACNN and mACNN as coarse propagators can be understood through three complementary observations. First, the use of a larger time-step ratio $c$ is feasible because the neural coarse solver is not constrained by the stability restrictions that limit traditional numerical coarse solvers. The ACNN and mACNN models are trained with data generated at a temporal resolution, ensuring that the local dynamics are accurately represented. Within the Parareal framework, however, the network can be applied recursively to advance the solution over a larger effective time step $\Delta t$ by composing multiple evaluations of the trained model. In this way, the neural coarse solver remains stable and accurate while enabling significantly more efficient coarse propagation than a numerical solver, whose stability often deteriorates when large time steps are used. Second, the number of iterations $k$ required for convergence is directly linked to the accuracy of the coarse solver: since ACNN and mACNN can achieve relative higher accuracy, particularly in capturing sharp interfacial dynamics, fewer correction steps are needed, which further enhances efficiency. Third, the ratio $\tfrac{t_{\text{num}}}{t_{\text{nn}}}$ reflects the computational advantage of the neural coarse solver, whose evaluation cost is much faster compared to the fine solver. Beyond efficiency, these models also offer strong generalization capabilities: they can be applied to arbitrary initial conditions without retraining and can naturally handle nonlocal terms in the mAC equation. Together, these properties ensure that ACNN and mACNN provide an accurate, efficient, and broadly applicable coarse propagator for the Parareal algorithm.

\section{Numerical experiments}\label{sec:results}
In this section, we first illustrate a series of experiments  classic AC equation \eqref{eq:class_AC_eq} and mass-conserve AC equation \eqref{eq:mass_AC_eq} in 2D and 3D on the convergence and accuracy of Parareal algorithm. And then we perform the speedup on multiply GPUs comparing with fine solver time. Throughout all experiments, we fix the spatial domain $\Omega=[-0.5,0.5]^d$ with $d=2,3$, the inerfacial parameter $\epsilon=0.01$ for 2D examples and $\epsilon=0.02$ for 3D examples, and use $\Delta t =0.1$ for coarse function and $\Delta t^* =0.001$ for fine solver. To evaluate the accuracy of our Parareal algorithm, we calculate three types of errors under the relative $L_2$-norm
\begin{equation*}
    \mathcal{E}_{\mathcal{NN}} = \frac{\lVert \mathbf{U}_{\mathcal{NN}}-\mathbf{U}_{ref}\rVert_{2}}{\lVert\mathbf{U}_{ref}\rVert_{2}},
\end{equation*}
\begin{equation*}
    \mathcal{E}_{Para} = \frac{\lVert \mathbf{U}_{Para}-\mathbf{U}_{ref}\rVert_{2}}{\lVert\mathbf{U}_{ref}\rVert_{2}},
\end{equation*}
\begin{equation*}
    \mathcal{E}_{Num} = \frac{\lVert \mathbf{U}_{Num}-\mathbf{U}_{ref}\rVert_{2}}{\lVert\mathbf{U}_{ref}\rVert_{2}},
\end{equation*}
where $\mathbf{U}_{\mathcal{NN}}$ is the predicted solution generated by CNN models, $\mathbf{U}_{Par}$ is the updated simulation solutions after the Parareal algorithm, $\mathbf{U}_{Num}$ is the fine numerical solutions with $\Delta t^*$ and $\mathbf{U}_{ref}$ is the reference solution. The reference solution is computed with the same spatial mesh and the sufficiently small time step size $\Delta t=0.0001$. 
We conduct the numerical experiments in PyTorch, and all experiments are run on Frontier with a AMD GPU card with 64GB of memory. 

\subsection{Training Stability Across Models}
\begin{figure}[!ht]
    \centering
    \includegraphics[width=0.47\linewidth]{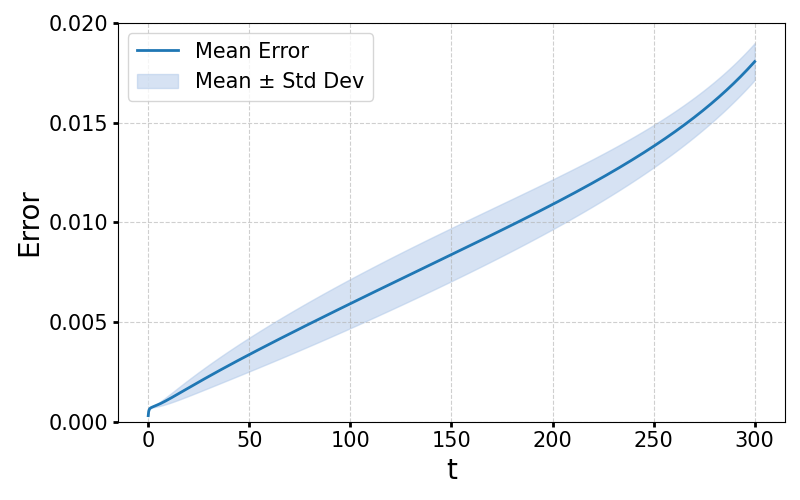}
    \includegraphics[width=0.47\linewidth]{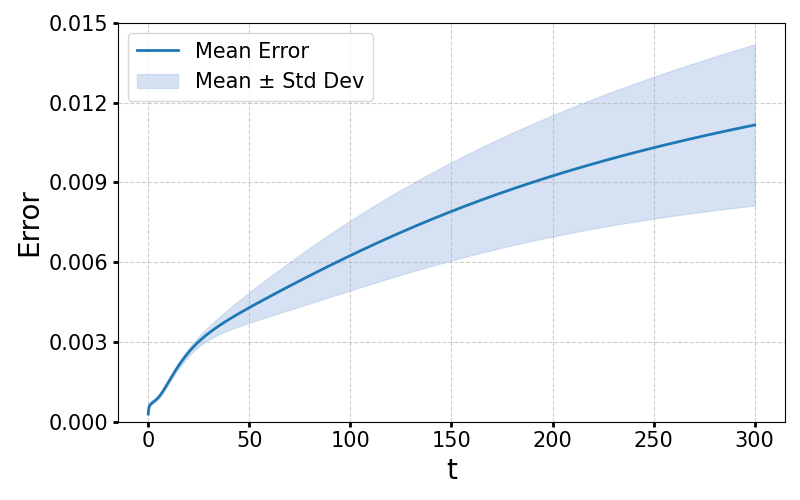}
    \caption{2D CNN models: Mean predicted errors and one-standard-deviation bands across ten independently trained models. The left panel shows results for ACNN, while the right panel presents results for mACNN.}
    \label{fig:2D_model_compare_errors}
\end{figure}

\begin{figure}[!ht]
    \centering
    \includegraphics[width=0.47\linewidth]{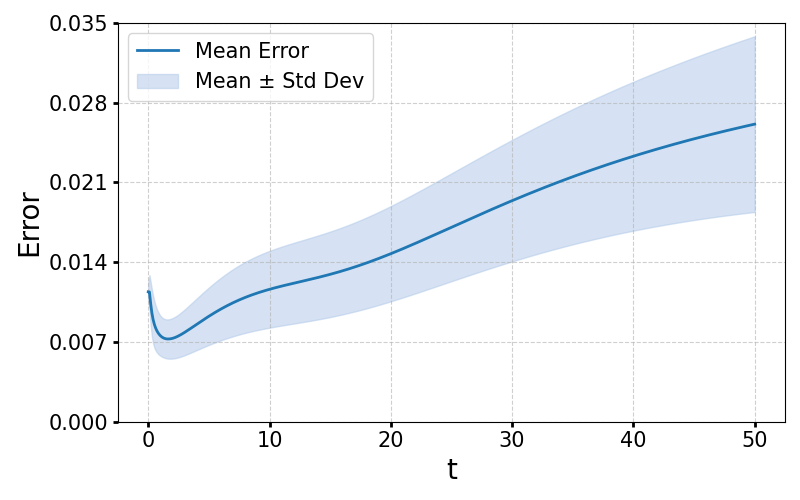}
    \includegraphics[width=0.47\linewidth]{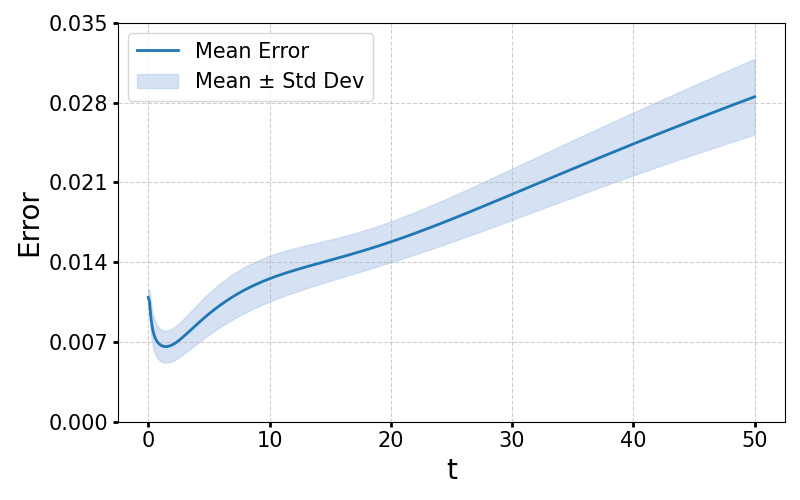}
    \caption{3D CNN models: Mean predicted errors and one-standard-deviation bands across ten independently trained models. The left panel corresponds to ACNN, and the right panel to mACNN.}
    \label{fig:3D_model_compare_errors}
\end{figure}

To systematically assess the robustness and training stability of our proposed models, we perform ten independent training runs for each architecture—ACNN and mACNN—under identical training data but with distinct random initializations. We then compute and visualize the mean predicted error along with the mean ± one standard deviation bands over time in both 2D (Figure \ref{fig:2D_model_compare_errors}) and 3D (Figure \ref{fig:3D_model_compare_errors}) settings. These results help quantify the variability due to stochastic training effects and offer insight into the consistency of the learned dynamics.

As illustrated in Figure \ref{fig:2D_model_compare_errors}, the 2D simulations show that both ACNN (left) and mACNN (right) produce highly stable and consistent predictions across all runs. The narrow error bands reflect low variance, indicating that the models converge reliably despite randomness in initialization or optimizer trajectories.

For the 3D case shown in Figure \ref{fig:3D_model_compare_errors}, the models exhibit slightly increased variability, as expected due to the higher dimensional complexity of the problem. However, even in this setting, the standard deviation remains modest, and the mean prediction error stays well below 0.035 throughout the entire time horizon. These results reinforce the reliability and generalization capabilities of our architectures.

Given the observed stability, we proceed with reporting detailed results based on a single representative model per architecture for the remainder of our experiments.

\subsection{Dynamics prediction for 2D examples}
Next, we present 2D benchmark examples to demonstrate the effectiveness of the Parareal algorithm in terms of both convergence and accuracy for the dynamics of the classical and conservative AC equations. Two well-known benchmark examples are considered: merging bubbles and grain coarsening.
For spatial discretization, we set the mesh size to $h=1/256$ with time step size $\Delta t=0.1$. 

\subsubsection{Merging bubbles}
\begin{figure}[!ht]
    \centering
        \includegraphics[width=0.48\textwidth]{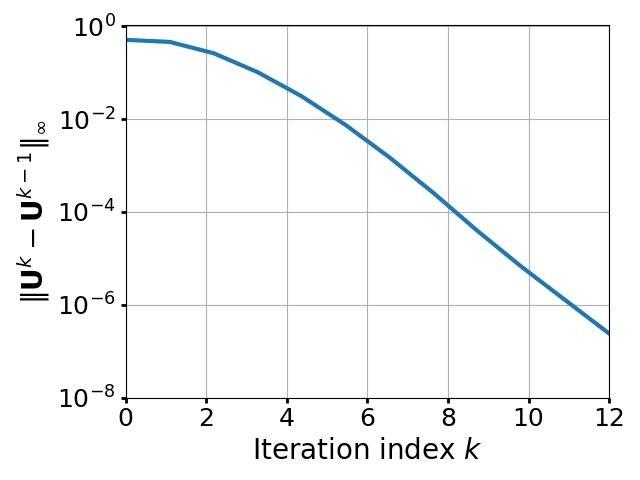}
        \includegraphics[width=0.475\textwidth]{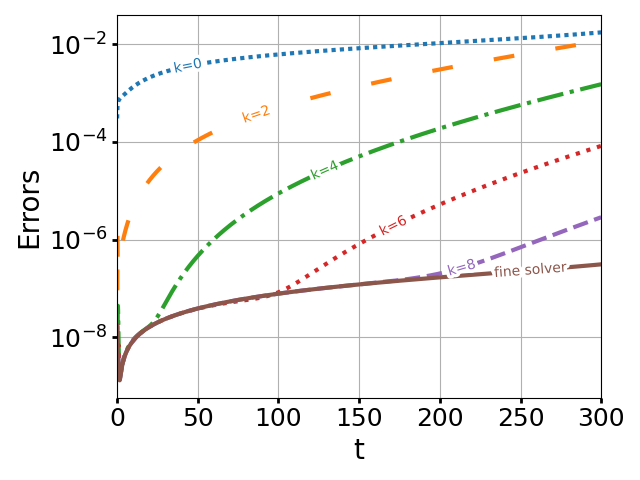}
    \caption{Classic AC in 2D: (Left) The $\|U^{k}-U^{k-1}\|_{\infty}$ with iteration index k and (right) the corresponding relative errors after the $k$-th iteration. Note $k=0$ means the coarse propagator errors without Parareal algorithm.}
    \label{fig:local_AC_2D_256}
\end{figure}

\begin{figure}[!ht]
    \centering
        \includegraphics[width=1\textwidth]{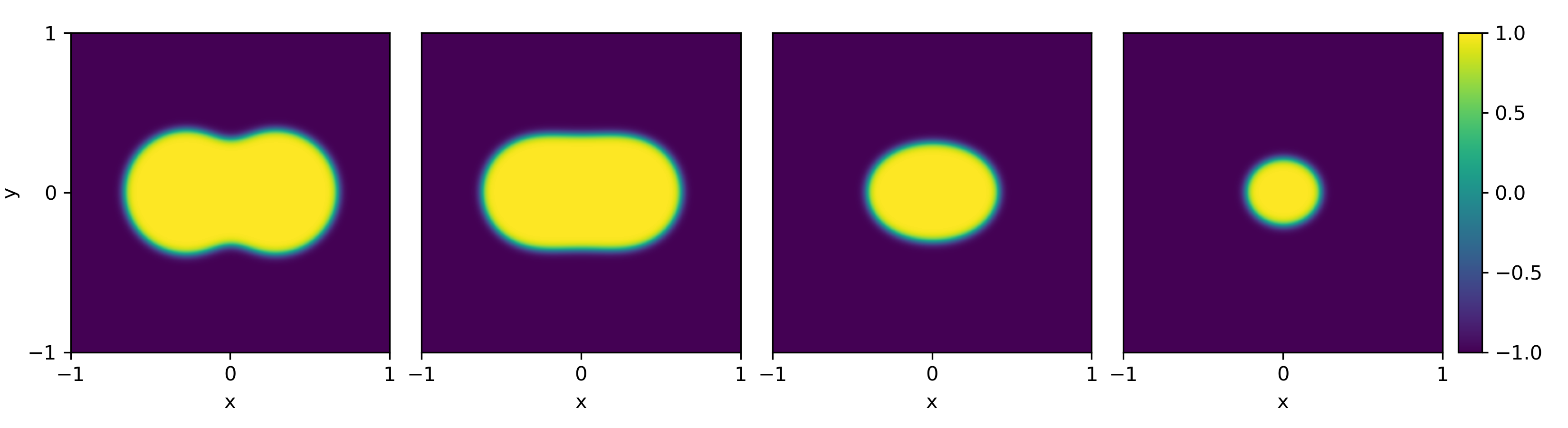}
        \includegraphics[width=1\textwidth]{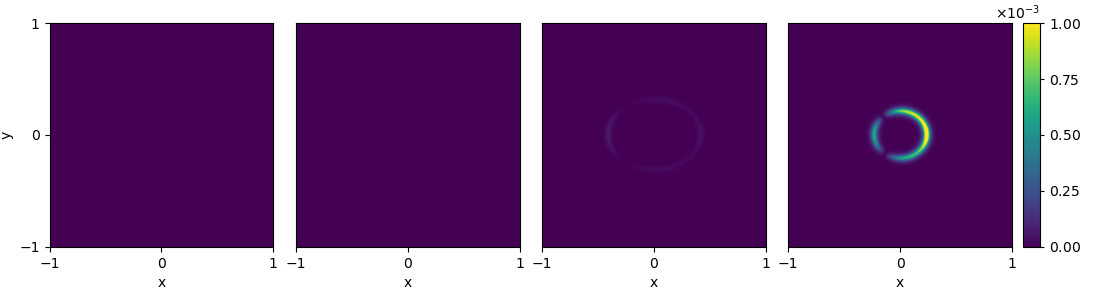}
    \caption{Classic AC  in 2D: The predicted solutions (top row) after six iterations of the Parareal algorithm and the corresponding numerical errors (bottom row) at $t=10,50,200,300$.}
    \label{fig:local_AC_2D_solutions}
\end{figure}

\begin{figure}[!ht]
    \centering
        \includegraphics[width=0.48\textwidth]{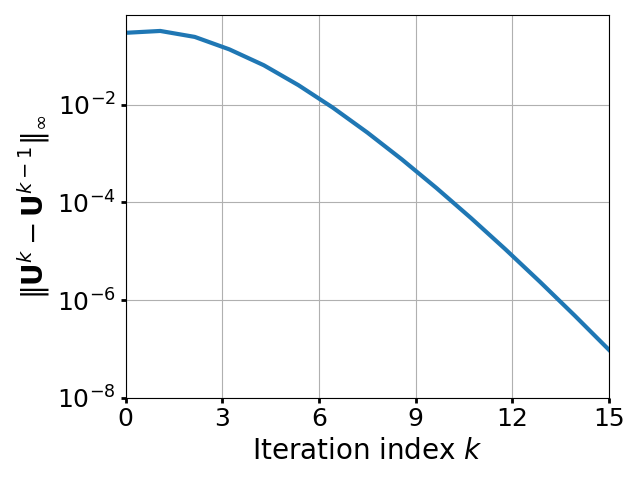}
        \includegraphics[width=0.475\textwidth]{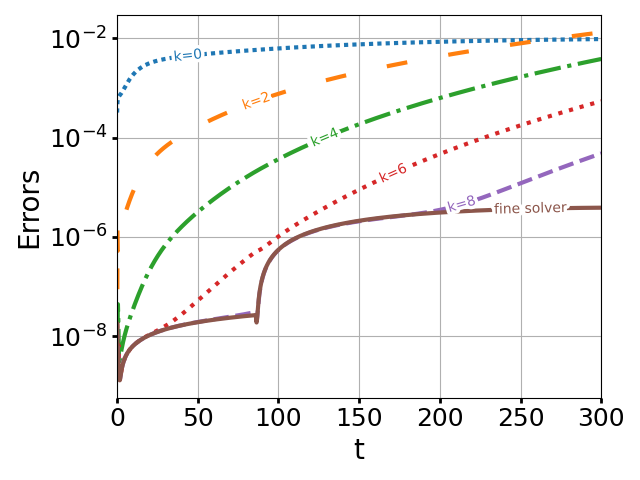}
    \caption{Mass-conserve AC  in 2D: (Left) The $\|U^{k}-U^{k-1}\|_{\infty}$ with iteration index k and (right) the corresponding relative errors after the $k$-th iteration. Note $k=0$ means the coarse propagator errors without Parareal algorithm.}
    \label{fig:local_AC_mass_2D_256}
\end{figure}

\begin{figure}[!ht]
    \centering
        \includegraphics[width=1\textwidth]{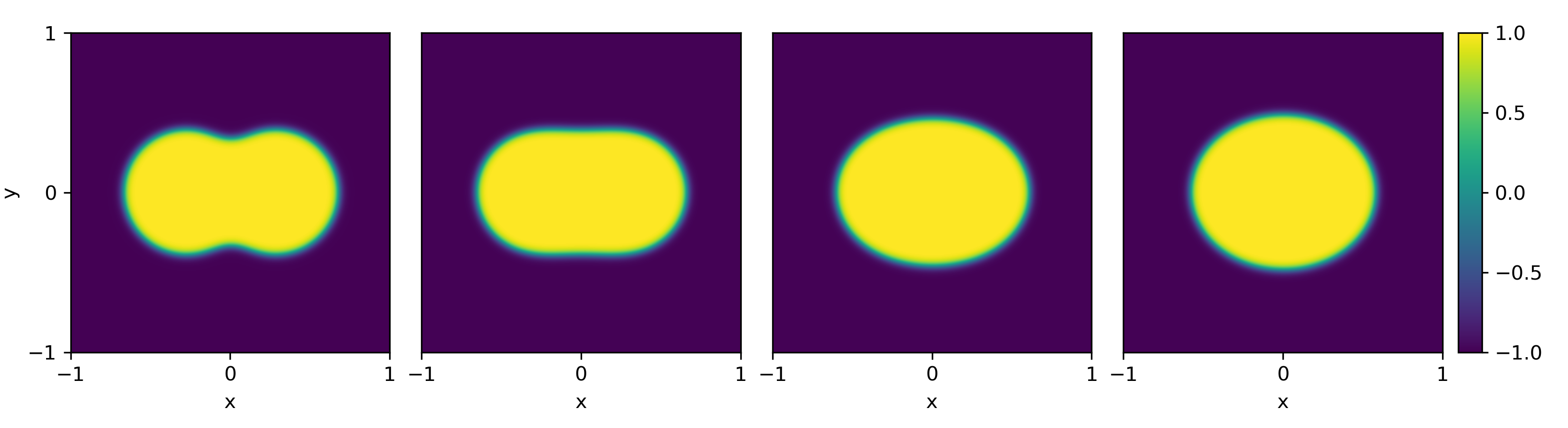}
        \includegraphics[width=1\textwidth]{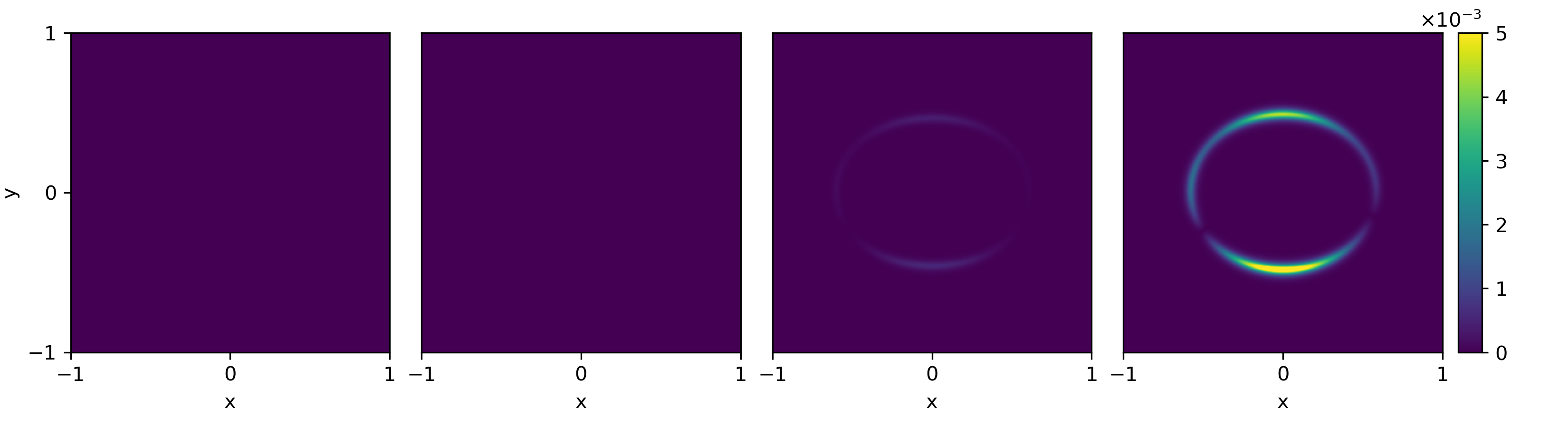}
    \caption{Mass-Conserving AC  in 2D: The predicted solutions (top row) after six iterations of the Parareal algorithm and the corresponding numerical errors (bottom row) at $t=10,50,200,300$.}
    \label{fig:local_AC_mass_2D_solutions}
\end{figure}

We consider the bubble merging problem, which is governed by the inital condition as follows: 
\begin{equation*}
    u_0(x,y)= \max\left(\tanh\left(\frac{0.2-\sqrt{(x-0.14)^2+y^2}}{\epsilon}\right), \tanh\left(\frac{0.2-\sqrt{(x+0.14)^2+y^2}}{\epsilon}\right) \right).
\end{equation*}
We predict the solution up to $t=300$, resulting in a total of $s=3000$ time steps.
Figure \ref{fig:local_AC_2D_256} and \ref{fig:local_AC_mass_2D_256} illustrates the convergence behavior of the Parareal algorithm with respect to the iteration index $k$, and the corresponding prediction errors after the $k$ iteration. The results show that the supremum norm difference $\|U^{k}-U^{k-1}\|_{\infty}$ decreases rapidly, reaching $10^{-6}$ around 10 iterations. 
With each Parareal iteration $k$, the solution becomes increasingly close to the fine solver. After a sufficient number of iterations (e.g., $k$=8), the error becomes negligible, showing convergence. 
Figure \ref{fig:local_AC_2D_solutions}, \ref{fig:local_AC_mass_2D_solutions} plots the predicted solution after six iterations of the Parareal algorithm and associated errors at the times $t=10,50,200,300$ for classic and mass-conservative AC equations, respectively.
The observed prediction errors are mainly localized in the interfacial regions between phases, and the maximum error compared to the reference is around 0.005.


\subsubsection{Grain coarsening}
\begin{figure}[!ht]
    \centering
        \includegraphics[width=0.48\textwidth]{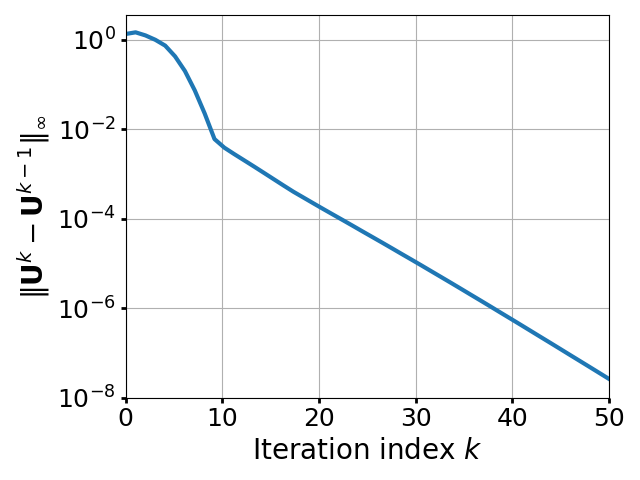}
        \includegraphics[width=0.475\textwidth]{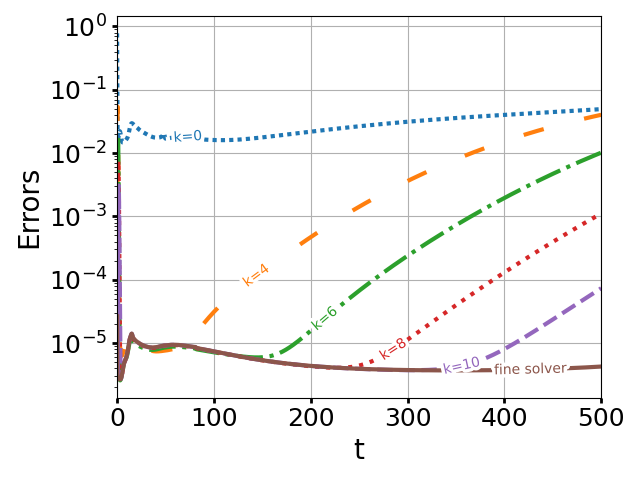}
        
    \caption{Classic AC  in 2D: (Left) The $\|U^{k}-U^{k-1}\|_{\infty}$ with iteration index k and (right) the corresponding relative errors after the $k$-th iteration. Note $k=0$ means the coarse propagtor errors without Parareal algorithm.}
    \label{fig:local_AC_2D_256_random}
\end{figure}

\begin{figure}[!ht]
    \centering
        \includegraphics[width=1\textwidth]{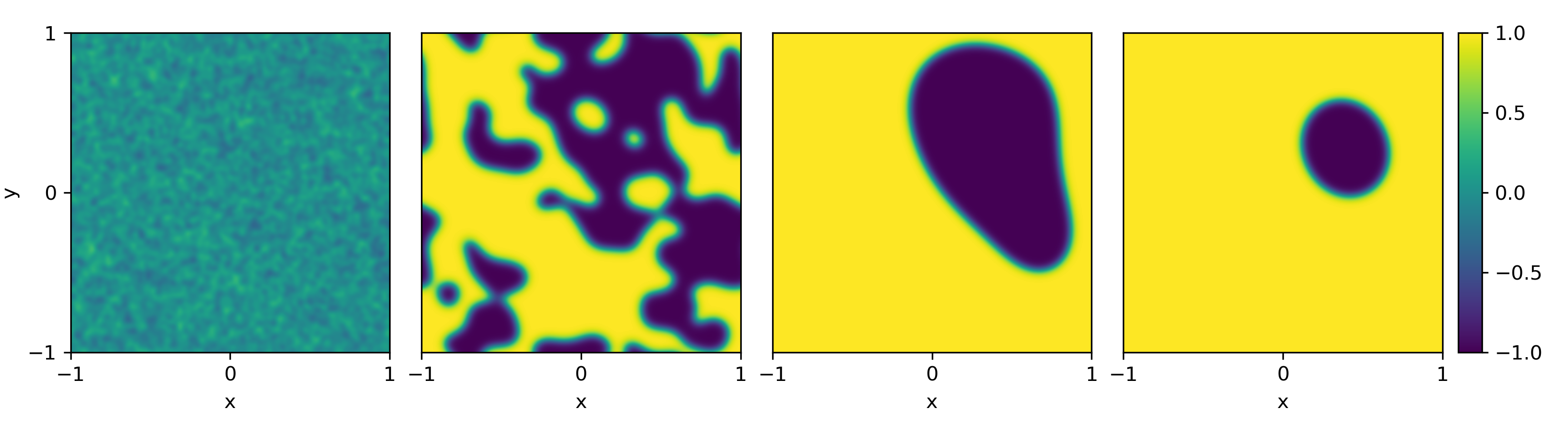}
        \includegraphics[width=1\textwidth]{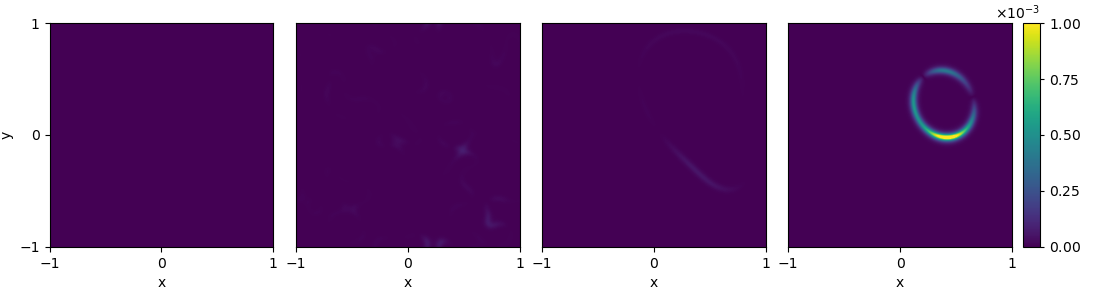}
    \caption{Classic AC  in 2D: The predicted solutions (top row) after six iterations of the Parareal algorithm and the corresponding numerical errors (bottom row) at $t=0.5,10,200,500$.}
    \label{fig:local_AC_2D_random_solutions}
\end{figure}

\begin{figure}[!ht]
    \centering
        \includegraphics[width=0.48\textwidth]{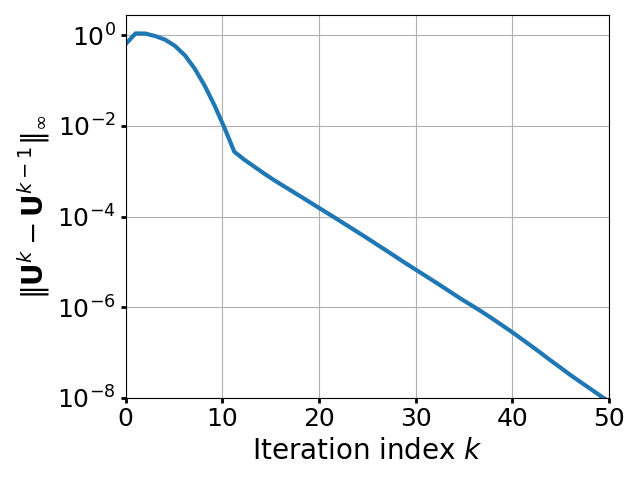}
        \includegraphics[width=0.475\textwidth]{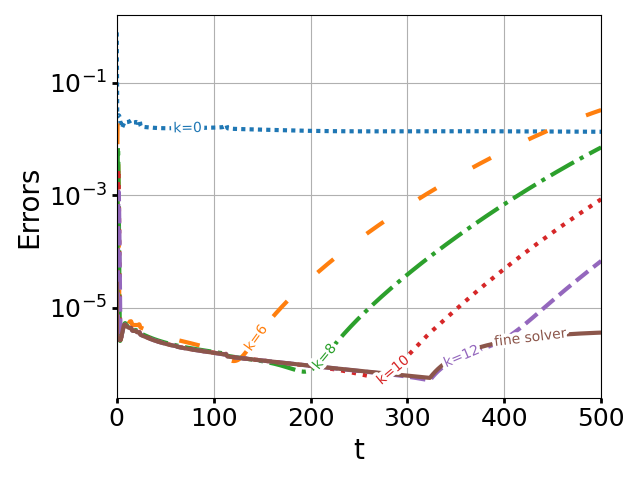}
    \caption{Mass-conserve AC  in 2D: (Left) The $\|U^{k}-U^{k-1}\|_{\infty}$ with iteration index k and (right) the corresponding relative errors after the $k$-th iteration. Note $k=0$ means the coarse propagator errors without Parareal algorithm.}
    \label{fig:local_AC_mass_2D_256_random}
\end{figure}

\begin{figure}[!ht]
    \centering
        \includegraphics[width=1\textwidth]{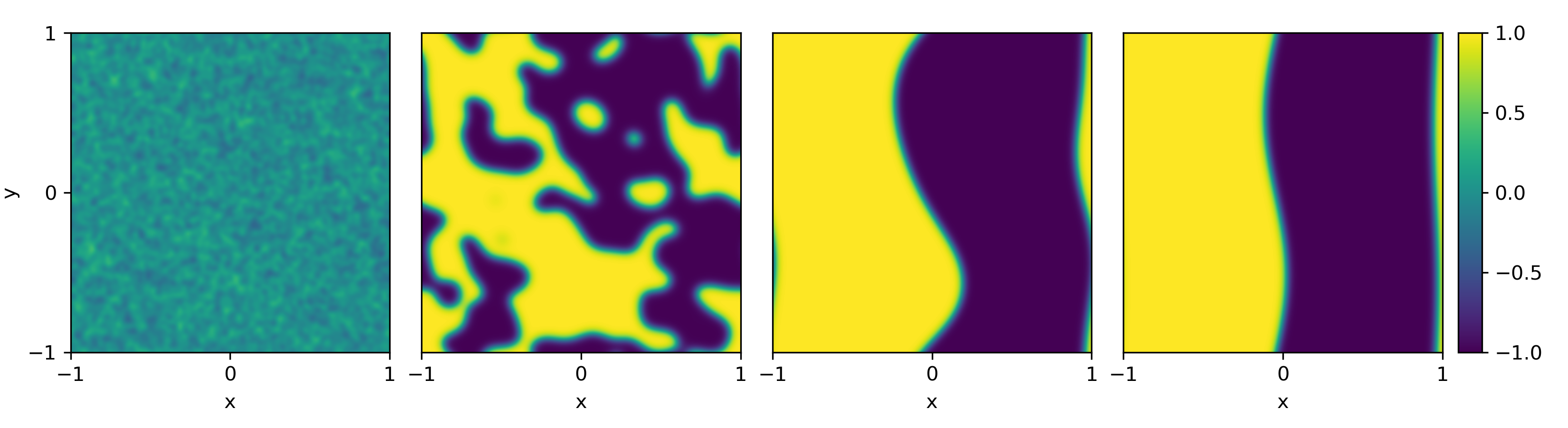}
        \includegraphics[width=1\textwidth]{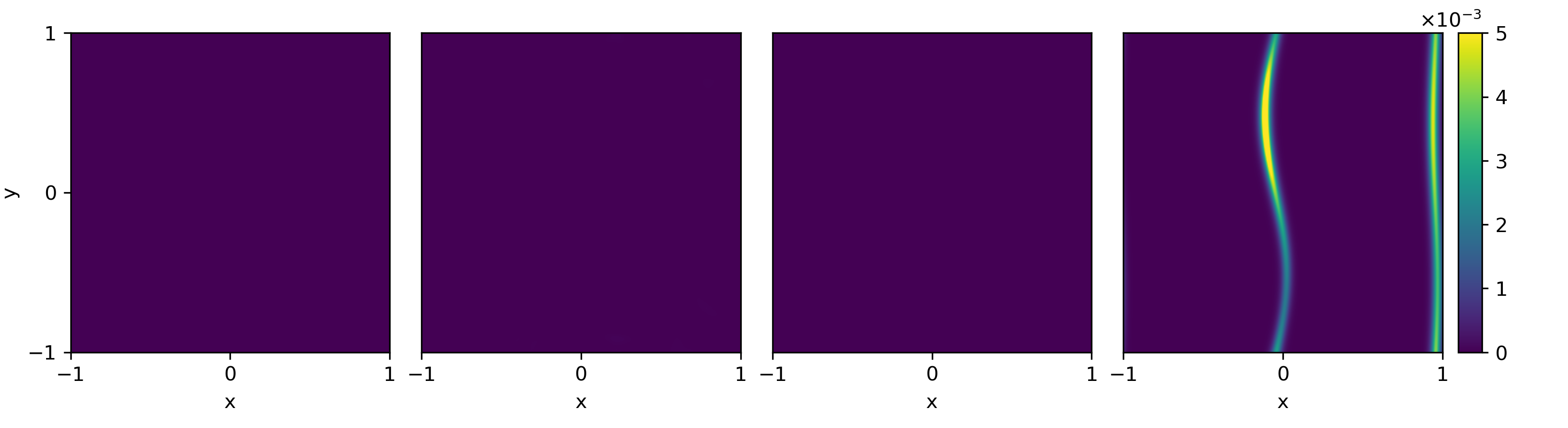}
    \caption{Mass-Conserving AC  in 2D: The predicted solutions (top row) after six iterations of the Parareal algorithm and the corresponding numerical errors (bottom row) at $t=0.5,10,200,500$.}
    \label{fig:local_AC_mass_2D_random_solutions}
\end{figure}

Next, we consider a random initialized condition of the form $u_0(x,y) = 0.9 \textbf{rand}(\cdot)$, and test it on classical and mass-conservative AC equations. For the case of grain coarsening, we predict the solution up to $t=500$, resulting in a total of $s=5000$ time steps. Figures~\ref{fig:local_AC_2D_256_random} and~\ref{fig:local_AC_mass_2D_256_random} (left) show the convergence behavior for both cases. Compared to the merging bubbles example, the convergence is slower, however, the same accuracy trend is observed in Figures~\ref{fig:local_AC_2D_256_random} and~\ref{fig:local_AC_mass_2D_256_random} (right), where the predicted solutions achieve nearly the same accuracy as the fine solver around eight iterations. We plot in
Figures~\ref{fig:local_AC_2D_random_solutions} and~\ref{fig:local_AC_mass_2D_random_solutions} the Parareal-predicted solutions after six iterations and their corresponding error maps at times $t=0.5,10,200,500$ for the classical and mass-conservative AC equations, respectively. During the grain coarsening process, the classical model exhibits gradual shrinkage and eventual disappearance of one phase (the yellow region), while the mass-conservative model preserves the overall phase volume, with no phase completely vanishing. Compared with reference solutions, the maximum predicted error is around 0.005, and the errors occur only in the transition between the two phases.

\subsection{Dynamics prediction for 3D examples}
In this section, we present 3D benchmark examples to demonstrate the effectiveness of the Parareal algorithm in terms of both convergence and accuracy for the dynamics of AC equations. For spatial discretization, we set the mesh size to $h=1/128$ with time-step size $\Delta t =0.1$.
\subsubsection{Star shape}
\begin{figure}[!ht]
    \centering
        \includegraphics[width=0.48\textwidth]{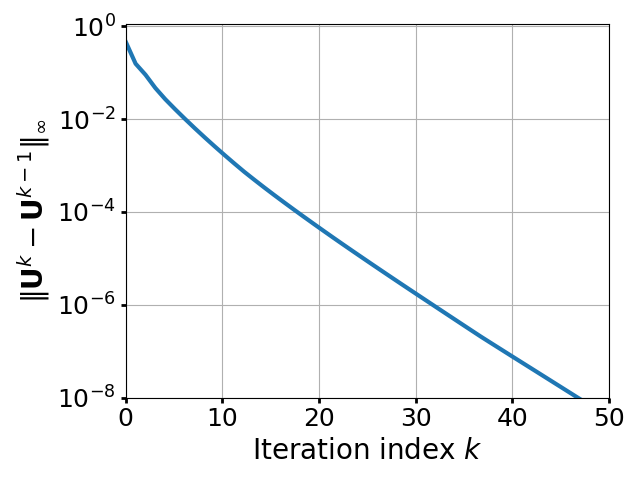}
        \includegraphics[width=0.475\textwidth]{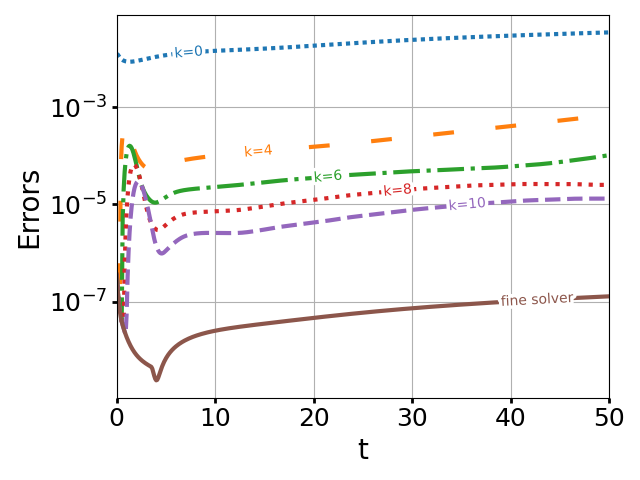}
    \caption{Classic AC  in 3D: (Left) The $\|U^{k}-U^{k-1}\|_{\infty}$ with iteration index k and (right) the corresponding relative errors after the $k$-th iteration. Note $k=0$ means the coarse propagator errors without the Parareal algorithm.}
    \label{fig:local_AC_3D}
\end{figure}

\begin{figure}[!ht]
   \centerline{\hspace{-0.55cm}
        \includegraphics[width=1\textwidth]{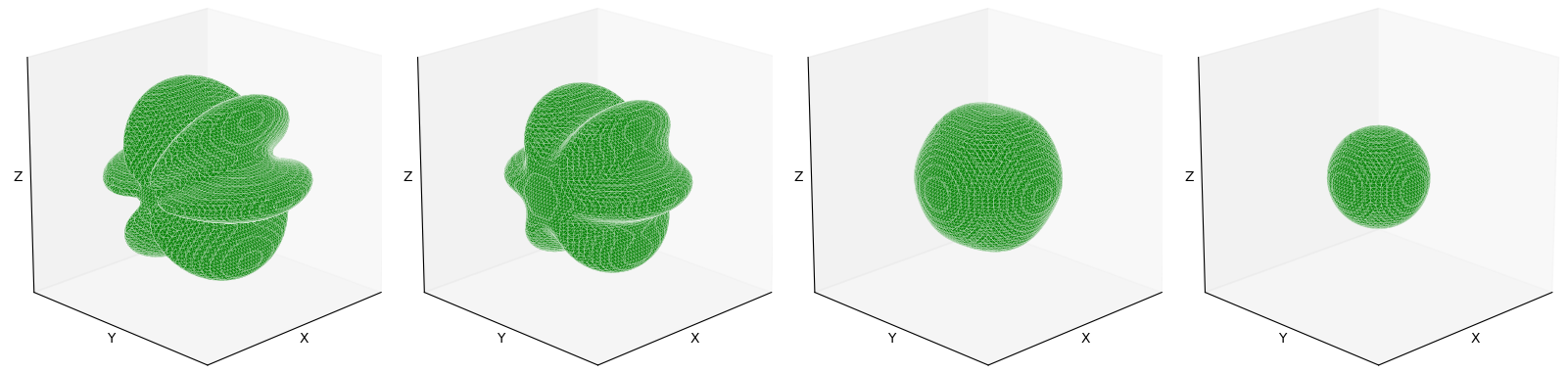}}
        \centerline{\hspace{-0.25cm}        \includegraphics[width=1\textwidth]{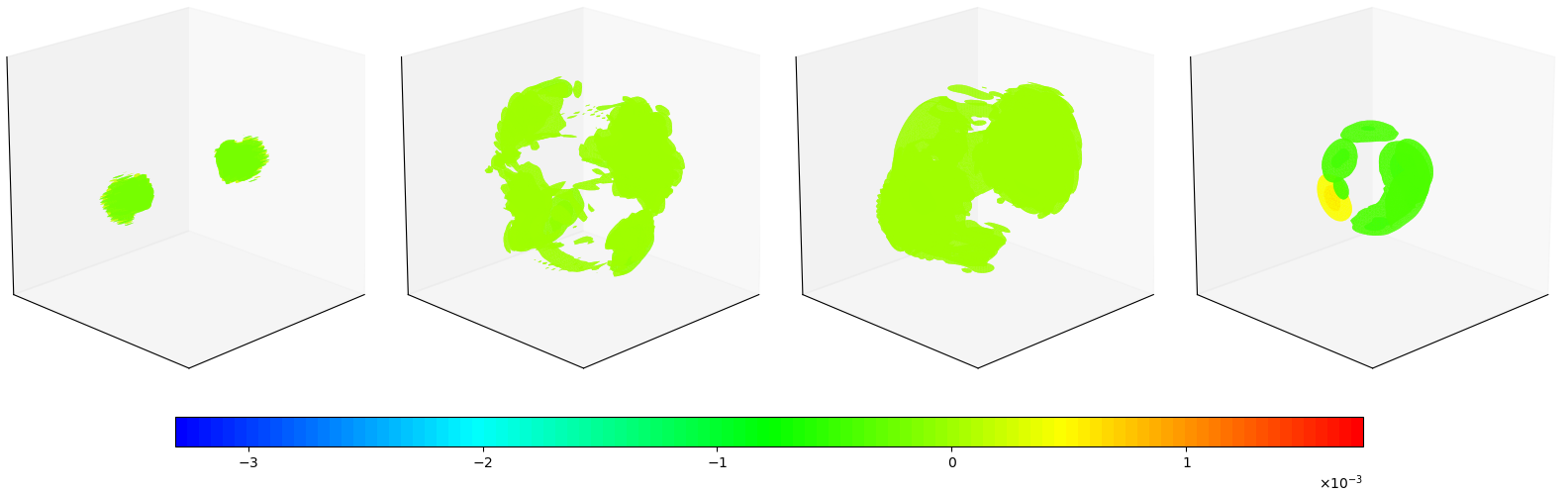}}
    \caption{Classic AC  in 3D: The predicted solutions (top row) after six iterations of the Parareal algorithm and the corresponding numerical errors (bottom row) at $t=1,10,20,50$.}
    \label{fig:local_AC_3D_solutions}
\end{figure}

\begin{figure}[!ht]
    \centering
        \includegraphics[width=0.48\textwidth]{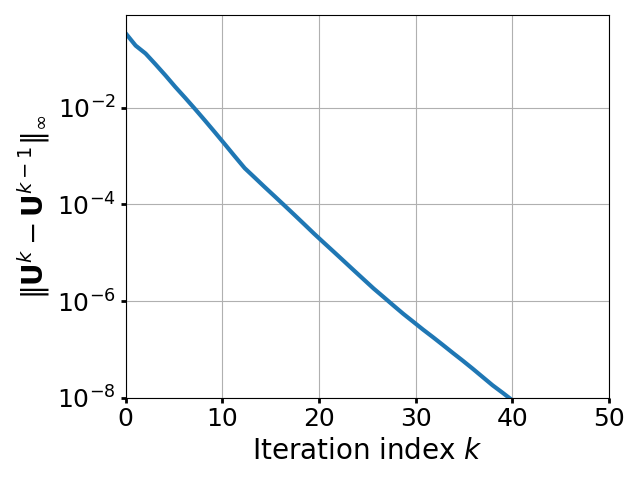}
        \includegraphics[width=0.475\textwidth]{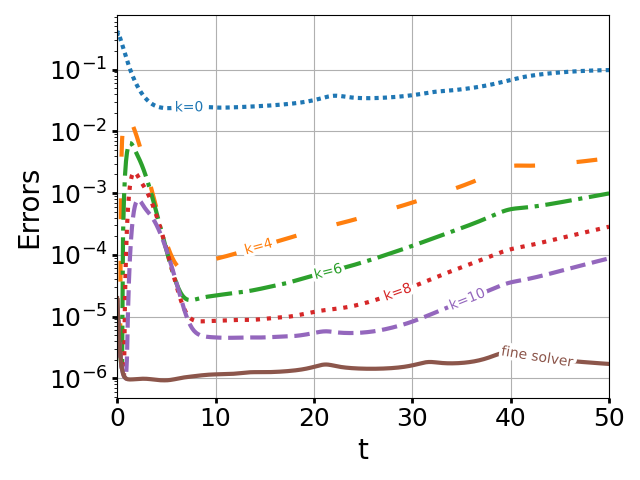}
    \caption{Classic AC  in 3D: (Left) The $\|U^{k}-U^{k-1}\|_{\infty}$ with iteration index k and (right) the corresponding relative errors after the $k$-th iteration. Note $k=0$ means the coarse propagator errors without the Parareal algorithm.}
    \label{fig:local_AC_mass_3D}
\end{figure}

\begin{figure}[!ht]
    \centerline{\hspace{-0.45cm}
        \includegraphics[width=1\textwidth]{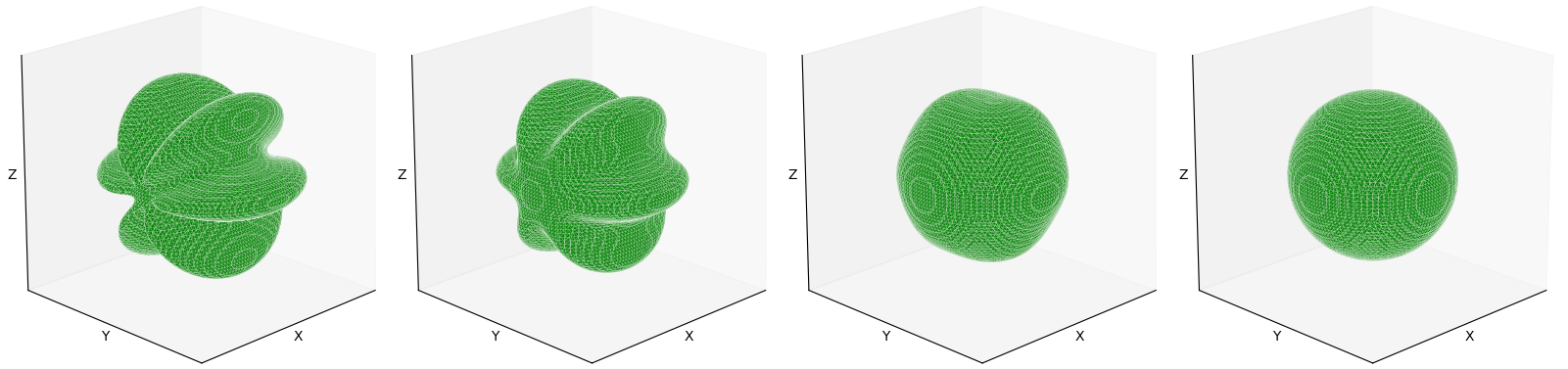}}
        \centerline{\hspace{-0.3cm} 
        \includegraphics[width=0.98\textwidth]{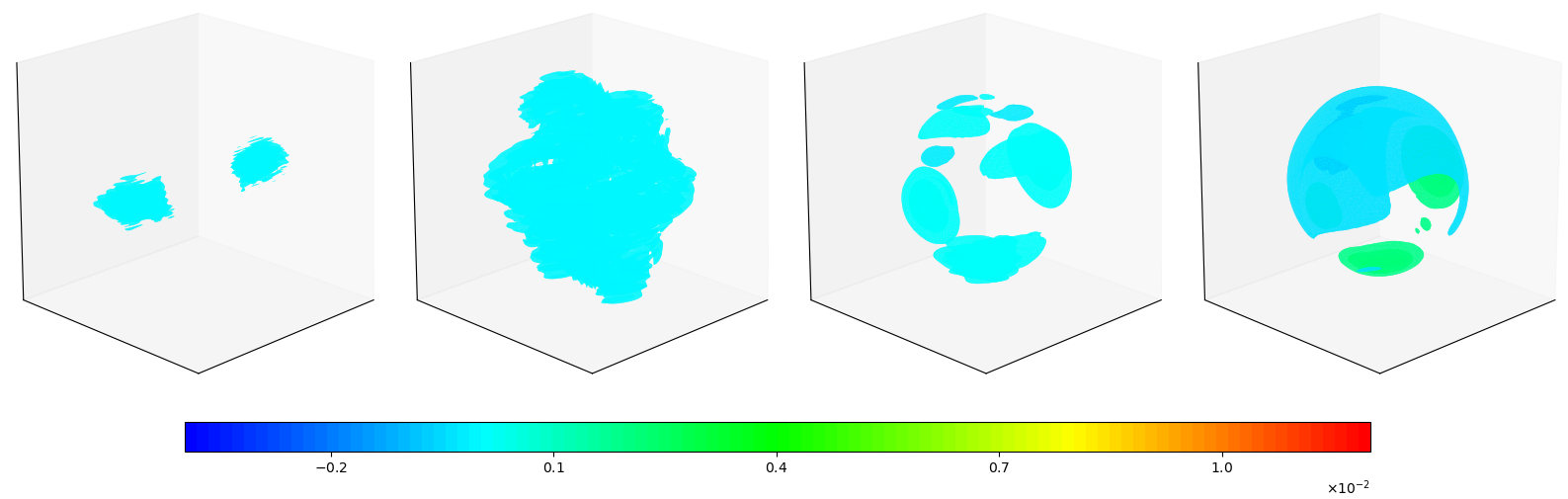}}
    \caption{Mass-Conserving AC  in 3D: The predicted solutions (top row) after six iterations of the Parareal algorithm and the corresponding numerical errors (bottom row) at $t=1,10,20,30$.}
    \label{fig:local_AC_mass_3D_solutions}
\end{figure}

We consider the "star" shape as the initial condition and predict the solution up to $t=50$ for the classical AC equation and up to $t=30$ for the mass-conservative AC equation, as the latter reaches a steady state by $t=30$.
Similarly to the 2D problem, we first test the convergence and accuracy of the Parareal algorithm on 3D models. Figures~\ref{fig:local_AC_3D} and~\ref{fig:local_AC_mass_3D} show the convergence and accuracy behavior for both cases. We observe that the convergence is fast, with the supremum norm difference $\|U^{k}-U^{k-1}\|_{\infty}$ reaching $10^{-8}$ after approximately 50 iterations. As the number of iterations increases, the accuracy improves and gradually approaches the solution obtained by the fine solver. Although convergence is slower than in the 2D cases, the predicted error reduces to around $10^{-5}$ after six iterations.
Figures \ref{fig:local_AC_3D_solutions} \ref{fig:local_AC_mass_3D_solutions} illustrate
The predicted solution after six iterations of the Parareal algorithm and the associated errors at times $t=1,10,20,50$ for the classic and $t=1,10,20,30$ mass-conservative AC equations, respectively. For the classic AC model, we observe that the initial star-shaped structure gradually smooths into a ball, which then progressively shrinks over time. In contrast, under the conservative AC setting, the ball eventually evolves into a stable, well-shaped structure with the same volume as the initial state, reflecting the mass conservation property inherent in the conservative formulation.
The maximum predicted error is around 0.01 compared to the reference solutions.



\subsubsection{Grain Coarsening}
\begin{figure}[!ht]
    \centering
        \includegraphics[width=0.48\textwidth]{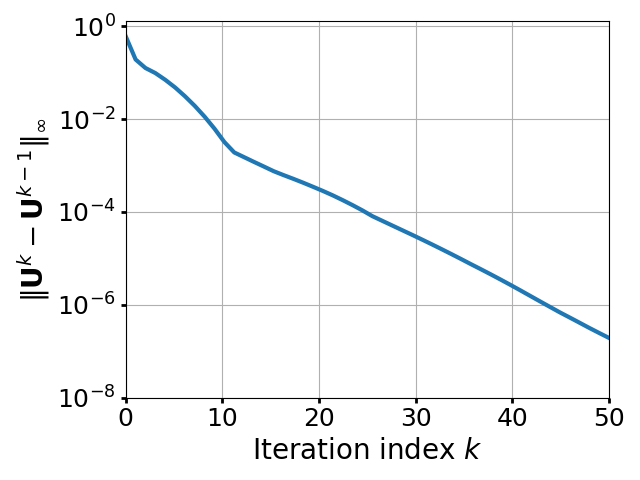}
        \includegraphics[width=0.475\textwidth]{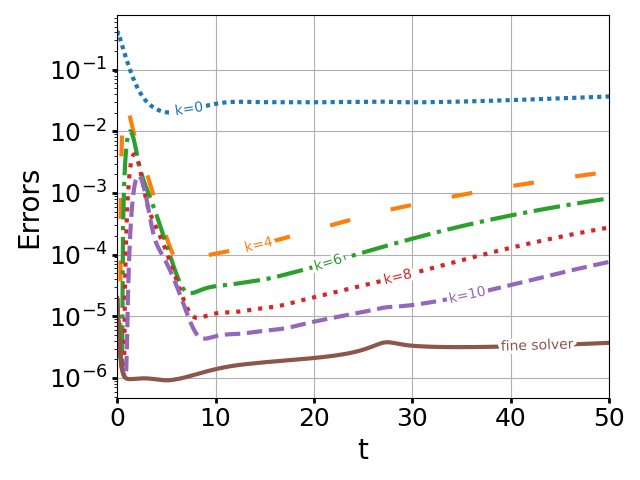}
    \caption{Classic AC in 3D: (Left) The $\|U^{k}-U^{k-1}\|_{\infty}$ with iteration index k and (right) the corresponding relative errors after the $k$-th iteration. Note $k=0$ means the coarse propagator errors without the Parareal algorithm.}
    \label{fig:random_local_AC_3D}
\end{figure}

\begin{figure}[!ht]
    \centerline{\hspace{-0.55cm}
        \includegraphics[width=1\textwidth]{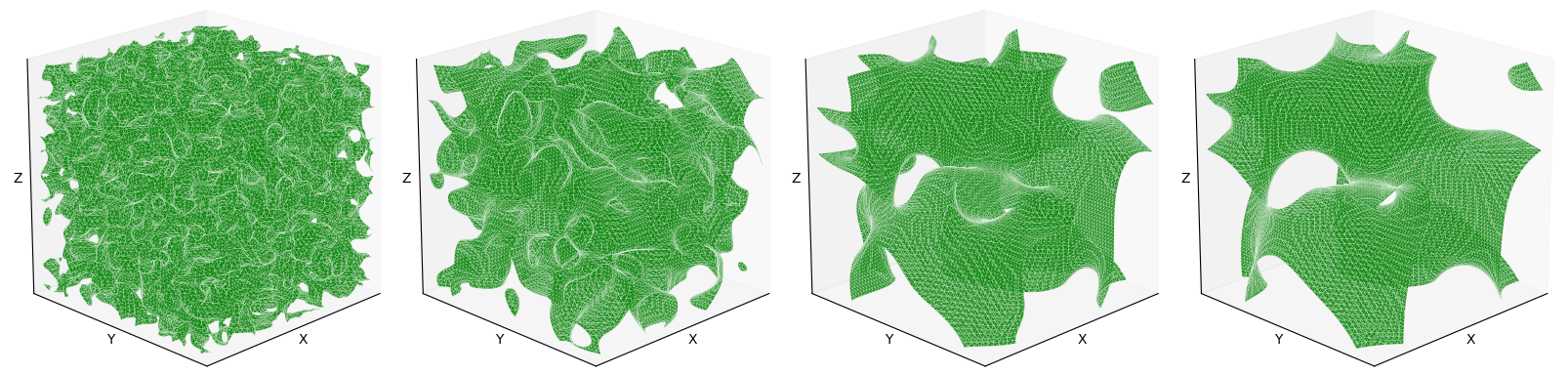}}
        \centerline{\hspace{-0.2cm} 
        \includegraphics[width=0.98\textwidth]{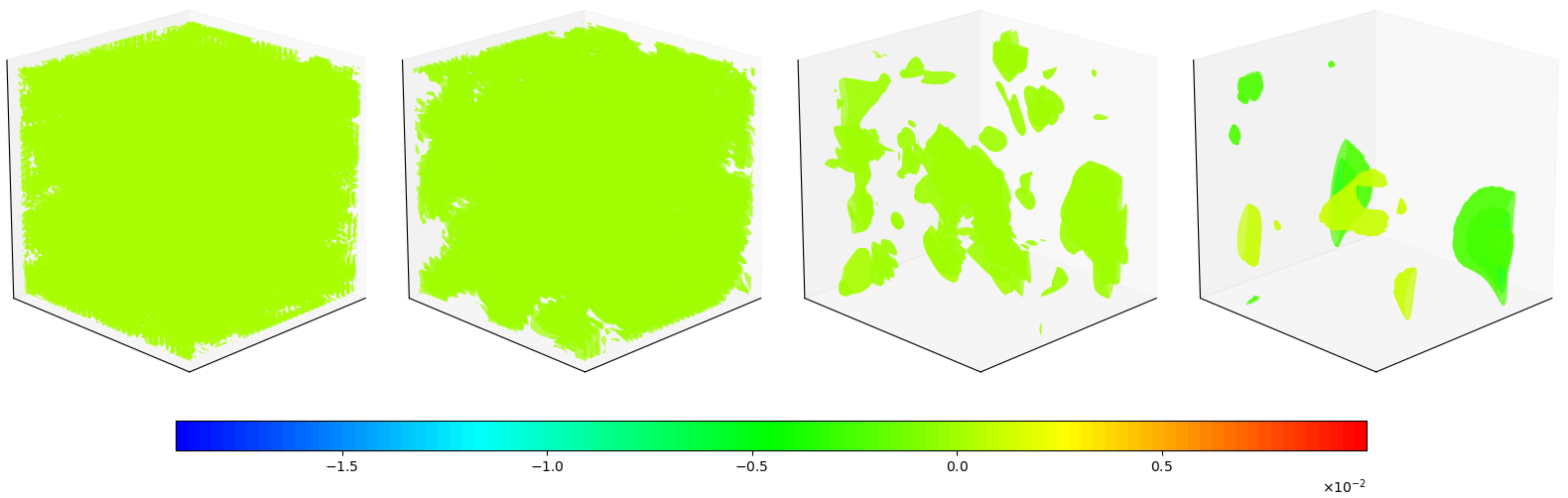}}
    \caption{Classic AC in 3D: The predicted solutions (top row) after six iterations of the Parareal algorithm and the corresponding numerical errors (bottom row) at $t=1,10,20,50$.}
    \label{fig:random_local_AC_3D_solutions}
\end{figure}

\begin{figure}[!ht]
    \centering
        \includegraphics[width=0.48\textwidth]{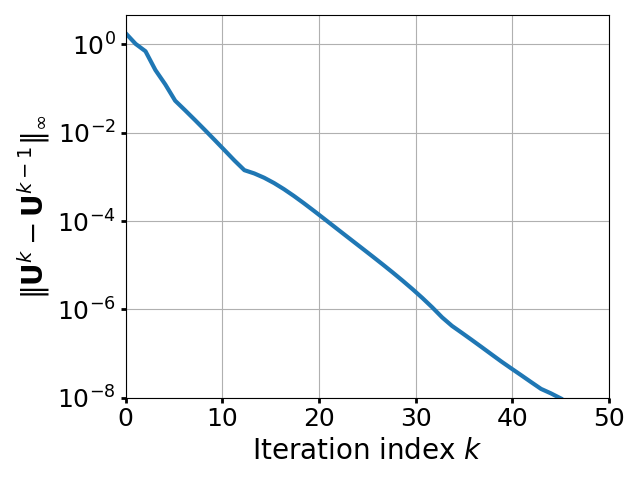}
        \includegraphics[width=0.475\textwidth]{figures_C4/random_local_AC_mass_3D_N64_T50_error.png}
    \caption{Classic AC in 3D: (Left) The $\|U^{k}-U^{k-1}\|_{\infty}$ with iteration index k and (right) the corresponding relative errors after the $k$-th iteration. Note $k=0$ means the coarse propagator errors without the Parareal algorithm.}
    \label{fig:random_local_AC_mass_3D}
\end{figure}

\begin{figure}[!ht]
    \centerline{\hspace{-0.6cm}
        \includegraphics[width=1\textwidth]{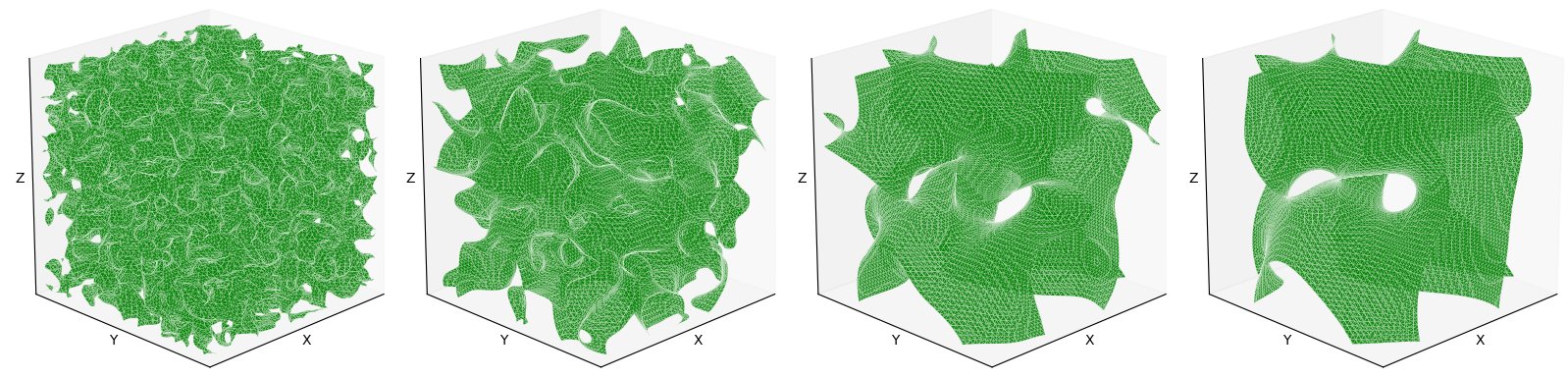}}
        \centerline{\hspace{-0.2cm}
        \includegraphics[width=0.98\textwidth]{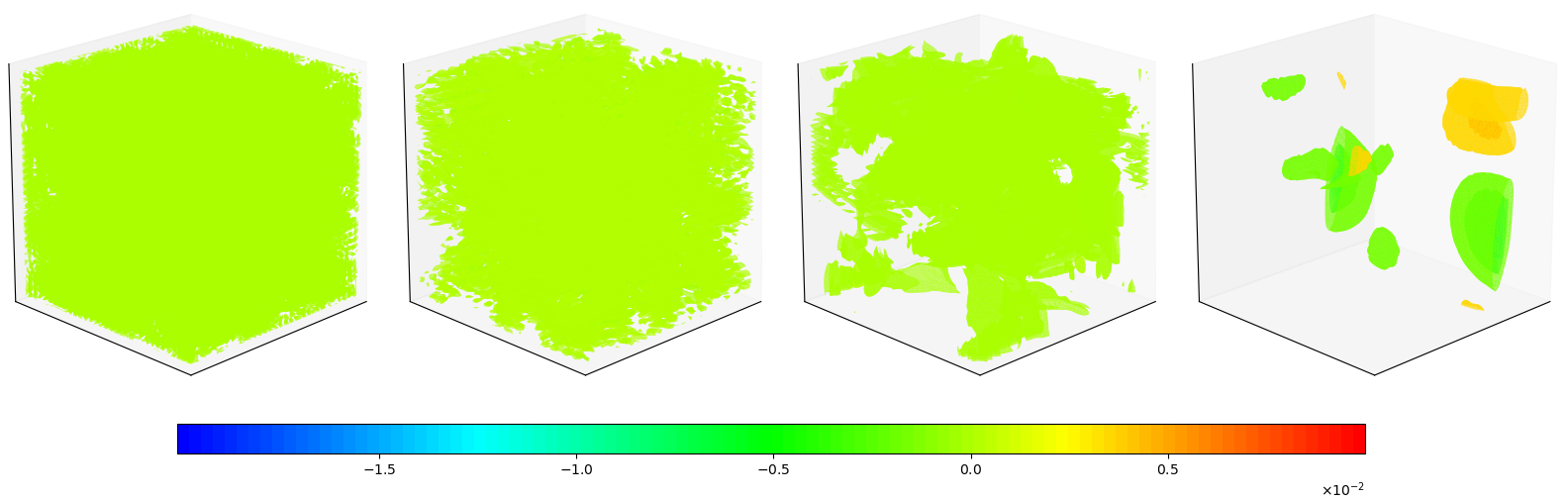}}
    \caption{Mass-Conserving AC in 3D: The predicted solutions (top row) after six iterations of the Parareal algorithm and the corresponding numerical errors (bottom row) at $t=1,10,20,50$.}
    \label{fig:random_local_AC_mass_3D_solutions}
\end{figure}

Finally, we apply the Parareal algorithm to predict the 3D grain coarsening dynamics with an intial condition generated randomly by $u_0(x,y,z) = 0.9\textbf{rand}(\cdot)$. The simulation time interval is set to be $[0,50]$ for both the classical and conservative AC equations. The convergence and precision with the iteration index $k$ are shown in Figures \ref{fig:random_local_AC_3D} and \ref{fig:random_local_AC_mass_3D}. The convergence and accuracy of the predicted solutions exhibit similar behavior to that observed in the previous example.
Figures \ref{fig:random_local_AC_3D_solutions} \ref{fig:random_local_AC_mass_3D_solutions} illustrate the predicted results with six iterations of the Parareal algorithm and associated errors at times $t=1,10,20,50$ for the classic and conservative AC equation, respectively.
Compared to the reference solution, the maximum errors are approximately 0.01.



\subsection{Speedup}
\begin{figure}[!ht]
    \centering
        \includegraphics[width=0.48\textwidth]{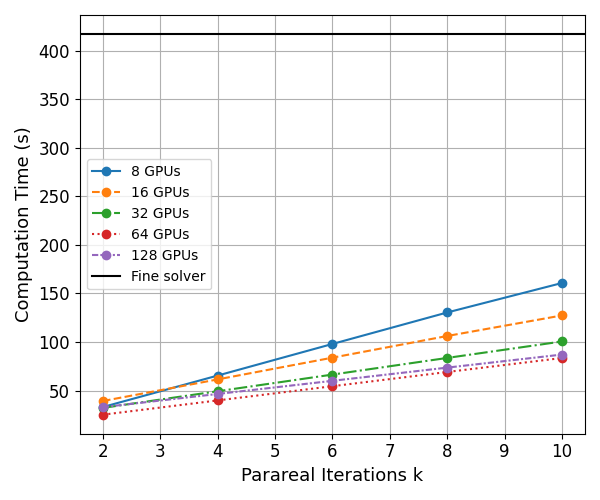}
        \includegraphics[width=0.48\textwidth]{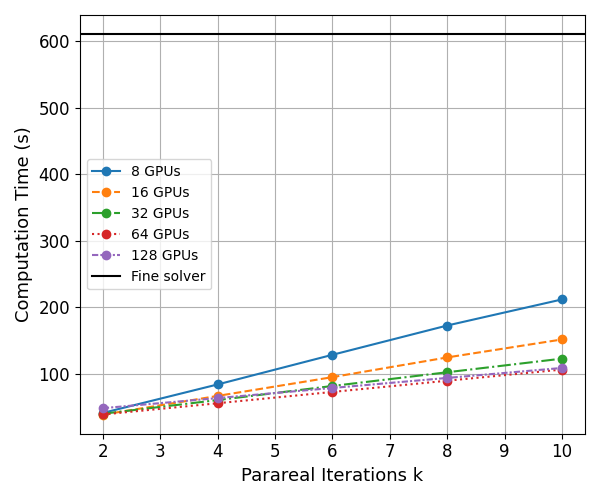}
    \caption{2D AC: The Parareal time with different iterations various on different number of GPUs for classica AC (left) and mass-conservative AC (right) equations.}
    \label{fig:2D_Parareal_time}
\end{figure}

\begin{figure}[!ht]
    \centering
        \includegraphics[width=0.48\textwidth]{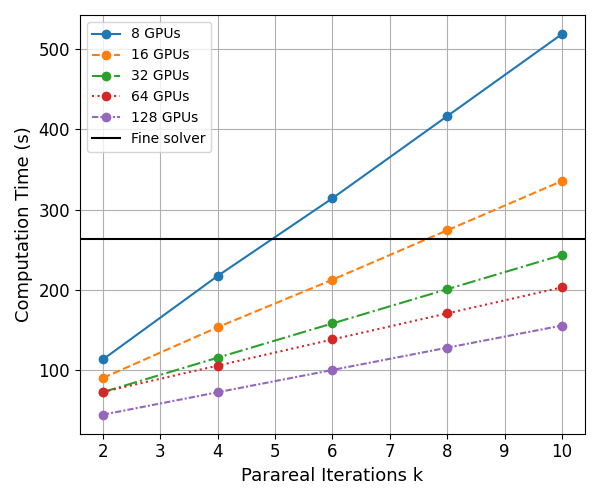}
        \includegraphics[width=0.48\textwidth]{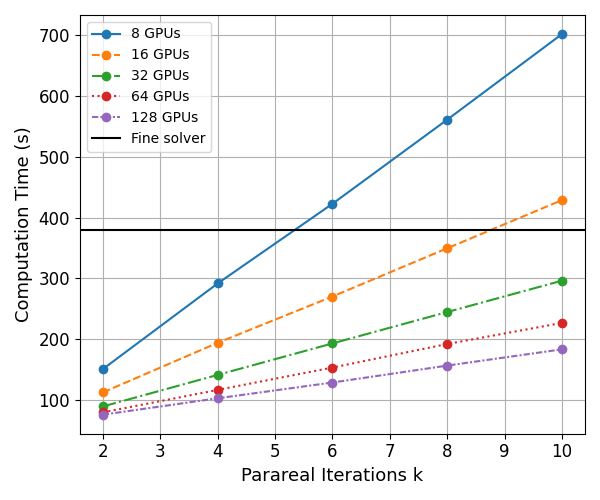}
    \caption{3D AC: The Parareal time with different iterations various on different number of GPUs for classica AC (left) and mass-conservative AC (right) equations.}
    \label{fig:3D_Parareal_time}
\end{figure}

In this part, we consider to estimate the speedup by applying the Parareal algorithm. To test that, we consider the random initial condition with ending time $T=500$ for 2D cases and ending time $T=50$ for 3D cases. Since the accuracy depends on the number of Parareal iterations, we present the computation times across multiple iteration counts and GPU counts in Figures~\ref{fig:2D_Parareal_time} and~\ref{fig:3D_Parareal_time}, and compare them to the computational cost of the fine solver. The different colored lines correspond to different numbers of GPUs (8, 16, 32, 64, and 128), and the black horizontal line indicates the computation time of the fine solver, used as a baseline reference. The goal of the Parareal algorithm is to significantly reduce simulation time while maintaining accuracy comparable to that of the fine solver. As expected, we observed that increasing the number of GPUs generally reduces the computation time.

For the 2D problems in Figure~\ref{fig:2D_Parareal_time}, as discussed in the previous section, the Parareal algorithm achieves high accuracy after eight iterations, closely matching the fine solver. At this point, the computation time is significantly reduced—by approximately a factor of six for both the classical and mass-conservative AC equations, when a sufficient number of GPUs is used. This demonstrates the algorithm’s effectiveness in accelerating time-dependent simulations without sacrificing accuracy. However, the benefit of adding more GPUs begins to diminish as a result of increasing communication overhead.

In the 3D experiments figure \ref{fig:3D_Parareal_time}, the Parareal algorithm still offers computational advantages, particularly when using a large number of GPUs (e.g., 64 or 128). However, the speedup is less pronounced compared to the 2D case, primarily because the simulation was only performed up to $T=50$ using 500 time steps, which limits the overall parallel time horizon. In such shorter time integration scenarios, the overhead from inter-GPU communication can outweigh the benefits of parallel execution. Even so, the Parareal method still achieves around a threefold speedup for both 3D variants of the Allen–Cahn equation when using sufficient GPUs.

Overall, the results show that with a well-trained coarse propagator and adequate hardware resources, the Parareal algorithm can achieve substantial acceleration in both 2D and 3D settings, while preserving the accuracy of the fine solver, even across various iteration counts and GPU configurations.


\section{Conclusions}
This work demonstrates the effectiveness of using a CNN as the coarse propagator in the Parareal algorithm for accelerating the simulation of both classical and mass-conservative AC equations. Numerical results in both 2D and 3D settings show that the proposed framework achieves significant speedup while maintaining accuracy comparable to the fine solver, especially when a sufficient number of GPUs is employed. The Parareal iterations converge rapidly, and in the 2D case, accuracy is attained within eight iterations, achieving up to 4–6 times runtime reduction. Although the speedup in 3D is less pronounced due to a shorter simulation horizon, the results still validate the scalability and robustness of the approach.
Looking forward, this method has strong potential for applications in complex multiscale systems such as weather forecasting, where a fast AI model can serve as the coarse propagator while the physical model—being too slow to run frequently—acts as the fine correction step. Unlike traditional Parareal settings that require a physics-based coarse solver, our approach removes this requirement by fully replacing the coarse phase solver with a trained neural network. This direction offers a promising path to real-time forecasting in scientific computing domains that currently face strict runtime limitations.

\section{Acknowledgment}
This work is supported by the U.S. Department of Energy, Office of Science, Office of Advanced Scientific
Computing Research, Applied Mathematics program, under the contracts ERKJ388 and ERKJ443.
ORNL is operated by UT-Battelle, LLC., for the U.S. Department of Energy under Contract DE-AC05-
00OR22725. Lili Ju acknowledges the support from the U.S. Department of Energy, Office of Science, Office of Advanced Scientific Computing Research program under grants DE-SC0025527 and DE-SC0022254. E. C. Cyr acknowledges support through funding from the U.S. Department of Energy, Office of Science, Office of Advanced Scientific
Computing Research at Sandia National Laboratoriers.
Sandia National Laboratories is a multimission laboratory managed and operated by National Technology \& Engineering Solutions of Sandia, LLC, a wholly owned subsidiary of Honeywell International Inc., for the U.S. Department of Energy’s National Nuclear Security Administration under contract DE-NA0003525. This paper describes objective technical results and analysis. Any subjective views or opinions that might be expressed in the paper do not necessarily represent the views of the U.S. Department of Energy or the United States Government.

\section*{Data availibility}
Data available on request from the authors.
%
\section*{Declarations}
This manuscript has been authored by UT-Battelle, LLC, under contract DE-AC05-00OR22725 with the US Department of Energy (DOE). The US government retains and the publisher, by accepting the article for publication, acknowledges that the US government retains a nonexclusive, paid-up, irrevocable, worldwide license to publish or reproduce the published form of this manuscript, or allow others to do so, for US government purposes. DOE will provide public access to these results of federally sponsored research in accordance with the DOE Public Access Plan.

\bibliographystyle{acm}
\bibliography{ref}

@inproceedings{ibrahim2023parareal,
  title={Parareal with a physics-informed neural network as coarse propagator},
  author={Ibrahim, Abdul Qadir and G{\"o}tschel, Sebastian and Ruprecht, Daniel},
  booktitle={European Conference on Parallel Processing},
  pages={649--663},
  year={2023},
  organization={Springer}
}

@article{krishnapriyan2021characterizing,
  title={Characterizing possible failure modes in physics-informed neural networks},
  author={Krishnapriyan, Aditi and Gholami, Amir and Zhe, Shandian and Kirby, Robert and Mahoney, Michael W},
  journal={Advances in neural information processing systems},
  volume={34},
  pages={26548--26560},
  year={2021}
}

@article{geng2024deep,
  title={A deep learning method for the dynamics of classic and conservative Allen-Cahn equations based on fully-discrete operators},
  author={Geng, Yuwei and Teng, Yuankai and Wang, Zhu and Ju, Lili},
  journal={Journal of Computational Physics},
  volume={496},
  pages={112589},
  year={2024},
  publisher={Elsevier}
}

@article{lions2001resolution,
  title={R{\'e}solution d'EDP par un sch{\'e}ma en temps {\guillemotleft}parar{\'e}el{\guillemotright}},
  author={Lions, Jacques-Louis and Maday, Yvon and Turinici, Gabriel},
  journal={Comptes Rendus de l'Acad{\'e}mie des Sciences-Series I-Mathematics},
  volume={332},
  number={7},
  pages={661--668},
  year={2001},
  publisher={Elsevier}
}

@article{du2021maximum,
  title={Maximum bound principles for a class of semilinear parabolic equations and exponential time-differencing schemes},
  author={Du, Qiang and Ju, Lili and Li, Xiao and Qiao, Zhonghua},
  journal={SIAM review},
  volume={63},
  number={2},
  pages={317--359},
  year={2021},
  publisher={SIAM}
}

@article{li2021unconditionally,
  title={Unconditionally maximum bound principle preserving linear schemes for the conservative Allen--Cahn equation with nonlocal constraint},
  author={Li, Jingwei and Ju, Lili and Cai, Yongyong and Feng, Xinlong},
  journal={Journal of Scientific Computing},
  volume={87},
  pages={1--32},
  year={2021},
  publisher={Springer}
}

@article{rubinstein1992nonlocal,
  title={Nonlocal reaction—diffusion equations and nucleation},
  author={Rubinstein, Jacob and Sternberg, Peter},
  journal={IMA Journal of Applied Mathematics},
  volume={48},
  number={3},
  pages={249--264},
  year={1992},
  publisher={Oxford University Press}
}

@article{minion2011hybrid,
  title={A hybrid parareal spectral deferred corrections method},
  author={Minion, Michael},
  journal={Communications in Applied Mathematics and Computational Science},
  volume={5},
  number={2},
  pages={265--301},
  year={2011},
  publisher={Mathematical Sciences Publishers}
}

@article{gunther2020layer,
  title={Layer-parallel training of deep residual neural networks},
  author={Gunther, Stefanie and Ruthotto, Lars and Schroder, Jacob B and Cyr, Eric C and Gauger, Nicolas R},
  journal={SIAM Journal on Mathematics of Data Science},
  volume={2},
  number={1},
  pages={1--23},
  year={2020},
  publisher={SIAM}
}

@article{falgout2014parallel,
  title={Parallel time integration with multigrid},
  author={Falgout, Robert D and Friedhoff, Stephanie and Kolev, Tz V and MacLachlan, Scott P and Schroder, Jacob B},
  journal={SIAM Journal on Scientific Computing},
  volume={36},
  number={6},
  pages={C635--C661},
  year={2014},
  publisher={SIAM}
}

@inproceedings{kirby2020layer,
  title={Layer-parallel training with gpu concurrency of deep residual neural networks via nonlinear multigrid},
  author={Kirby, Andrew and Samsi, Siddharth and Jones, Michael and Reuther, Albert and Kepner, Jeremy and Gadepally, Vijay},
  booktitle={2020 IEEE high performance extreme computing conference (HPEC)},
  pages={1--7},
  year={2020},
  organization={IEEE}
}

@article{meng2020ppinn,
  title={PPINN: Parareal physics-informed neural network for time-dependent PDEs},
  author={Meng, Xuhui and Li, Zhen and Zhang, Dongkun and Karniadakis, George Em},
  journal={Computer Methods in Applied Mechanics and Engineering},
  volume={370},
  pages={113250},
  year={2020},
  publisher={Elsevier}
}

@article{lorin2020derivation,
  title={Derivation and analysis of parallel-in-time neural ordinary differential equations},
  author={Lorin, Emmanuel},
  journal={Annals of Mathematics and Artificial Intelligence},
  volume={88},
  number={10},
  pages={1035--1059},
  year={2020},
  publisher={Springer}
}

@article{yalla2018parallel,
  title={Parallel in time algorithms for multiscale dynamical systems using interpolation and neural networks},
  author={Yalla, Gopal R and Engquist, Bj{\"o}rn},
  journal={Institute for Computational Engineering and Sciences},
  year={2018}
}

@article{agboh2020parareal,
  title={Parareal with a learned coarse model for robotic manipulation},
  author={Agboh, Wisdom and Grainger, Oliver and Ruprecht, Daniel and Dogar, Mehmet},
  journal={Computing and visualization in science},
  volume={23},
  number={1},
  pages={8},
  year={2020},
  publisher={Springer}
}

@article{nguyen2023numerical,
  title={Numerical wave propagation aided by deep learning},
  author={Nguyen, Hieu and Tsai, Richard},
  journal={Journal of Computational Physics},
  volume={475},
  pages={111828},
  year={2023},
  publisher={Elsevier}
}

@article{feng2003numerical,
  title={Numerical analysis of the Allen-Cahn equation and approximation for mean curvature flows},
  author={Feng, Xiaobing and Prohl, Andreas},
  journal={Numerische Mathematik},
  volume={94},
  pages={33--65},
  year={2003},
  publisher={Springer}
}

@article{xu2019stability,
  title={On the stability and accuracy of partially and fully implicit schemes for phase field modeling},
  author={Xu, Jinchao and Li, Yukun and Wu, Shuonan and Bousquet, Arthur},
  journal={Computer Methods in Applied Mechanics and Engineering},
  volume={345},
  pages={826--853},
  year={2019},
  publisher={Elsevier}
}

@article{shen2010numerical,
  title={Numerical approximations of allen-cahn and cahn-hilliard equations},
  author={Shen, Jie and Yang, Xiaofeng},
  journal={Discrete Contin. Dyn. Syst},
  volume={28},
  number={4},
  pages={1669--1691},
  year={2010}
}

@article{tang2016implicit,
  title={Implicit-explicit scheme for the Allen-Cahn equation preserves the maximum principle},
  author={Tang, Tao and Yang, Jiang},
  journal={Journal of Computational Mathematics},
  pages={451--461},
  year={2016},
  publisher={JSTOR}
}

@article{allen1979microscopic,
  title={A microscopic theory for antiphase boundary motion and its application to antiphase domain coarsening},
  author={Allen, Samuel M and Cahn, John W},
  journal={Acta metallurgica},
  volume={27},
  number={6},
  pages={1085--1095},
  year={1979},
  publisher={Elsevier}
}

@article{wheeler1992phase,
  title={Phase-field model for isothermal phase transitions in binary alloys},
  author={Wheeler, Adam A and Boettinger, William J and McFadden, Geoffrey B},
  journal={Physical Review A},
  volume={45},
  number={10},
  pages={7424},
  year={1992},
  publisher={APS}
}

@article{benevs2004geometrical,
  title={Geometrical image segmentation by the Allen--Cahn equation},
  author={Bene{\v{s}}, Michal and Chalupeck{\`y}, Vladim{\i}́r and Mikula, Karol},
  journal={Applied Numerical Mathematics},
  volume={51},
  number={2-3},
  pages={187--205},
  year={2004},
  publisher={Elsevier}
}

@article{shao2010computational,
  title={Computational model for cell morphodynamics},
  author={Shao, Danying and Rappel, Wouter-Jan and Levine, Herbert},
  journal={Physical review letters},
  volume={105},
  number={10},
  pages={108104},
  year={2010},
  publisher={APS}
}

@article{yang2018uniform,
  title={UNIFORM L p-BOUND OF THE ALLEN-CAHN EQUATION AND ITS NUMERICAL DISCRETIZATION.},
  author={Yang, Jiang and Du, Qiang and Zhang, Wei},
  journal={International Journal of Numerical Analysis \& Modeling},
  volume={15},
  year={2018}
}

@article{targ2016resnet,
  title={Resnet in resnet: Generalizing residual architectures},
  author={Targ, Sasha and Almeida, Diogo and Lyman, Kevin},
  journal={arXiv preprint arXiv:1603.08029},
  year={2016}
}

@article{lee2022parareal,
  title={Parareal neural networks emulating a parallel-in-time algorithm},
  author={Lee, Youngkyu and Park, Jongho and Lee, Chang-Ock},
  journal={IEEE Transactions on Neural Networks and Learning Systems},
  year={2022},
  publisher={IEEE}
}

@article{mattey2022novel,
  title={A novel sequential method to train physics informed neural networks for Allen Cahn and Cahn Hilliard equations},
  author={Mattey, Revanth and Ghosh, Susanta},
  journal={Computer Methods in Applied Mechanics and Engineering},
  volume={390},
  pages={114474},
  year={2022},
  publisher={Elsevier}
}

@article{zhao2020solving,
  title={Solving Allen-Cahn and Cahn-Hilliard Equations using the Adaptive Physics Informed Neural Networks},
  author={Zhao, Colby L},
  journal={Communications in Computational Physics},
  volume={29},
  number={3},
  year={2020}
}

@article{wang2022modified,
  title={A modified physics informed neural networks for solving the partial differential equation with conservation laws},
  author={Wang, Haifeng and Qian, Xu and Sun, Yabing and Song, Songhe},
  journal={Available at SSRN 4274376},
  year={2022}
}

@article{chen2025learn,
  title={Learn sharp interface solution by homotopy dynamics},
  author={Chen, Chuqi and Yang, Yahong and Xiang, Yang and Hao, Wenrui},
  journal={arXiv e-prints},
  pages={arXiv--2502},
  year={2025}
}

@article{raissi2019physics,
  title={Physics-informed neural networks: A deep learning framework for solving forward and inverse problems involving nonlinear partial differential equations},
  author={Raissi, Maziar and Perdikaris, Paris and Karniadakis, George E},
  journal={Journal of Computational physics},
  volume={378},
  pages={686--707},
  year={2019},
  publisher={Elsevier}
}

@article{xu2025overview,
  title={Overview frequency principle/spectral bias in deep learning},
  author={Xu, Zhi-Qin John and Zhang, Yaoyu and Luo, Tao},
  journal={Communications on Applied Mathematics and Computation},
  volume={7},
  number={3},
  pages={827--864},
  year={2025},
  publisher={Springer}
}

@article{huang2025frequency,
  title={Frequency-adaptive multi-scale deep neural networks},
  author={Huang, Jizu and You, Rukang and Zhou, Tao},
  journal={Computer Methods in Applied Mechanics and Engineering},
  volume={437},
  pages={117751},
  year={2025},
  publisher={Elsevier}
}

@article{liu2024mitigating,
  title={Mitigating spectral bias for the multiscale operator learning},
  author={Liu, Xinliang and Xu, Bo and Cao, Shuhao and Zhang, Lei},
  journal={Journal of Computational Physics},
  volume={506},
  pages={112944},
  year={2024},
  publisher={Elsevier}
}

\end{document}